\documentclass[final,3p,number, sort]{elsarticle}

\usepackage{amssymb}
\usepackage{amsmath}
\usepackage{amsthm}
\usepackage{bm}
\usepackage[colorlinks=true]{hyperref}
\usepackage[section]{placeins}
\usepackage[T1]{fontenc}
\usepackage[bitstream-charter]{mathdesign} 
\usepackage[scaled]{berasans}              
\usepackage{color}
\usepackage{graphicx}
\usepackage[normalem]{ulem}

\biboptions{semicolon}

\newcommand{\AVE}[1]{\ensuremath{\langle {#1} \rangle}}
\newcommand{\ABS}[1]{\ensuremath{\lvert {#1} \rvert}}

\newcommand{\bI}{\ensuremath{\bm{I}}}

\newcommand{\bK}{\ensuremath{\bm{K}}}

\newcommand{\bP}{\ensuremath{\bm{P}}}

\newcommand{\bS}{\ensuremath{\bm{S}}}

\newcommand{\bU}{\ensuremath{\bm{U}}}
\newcommand{\bV}{\ensuremath{\bm{V}}}

\newcommand{\bff}{\ensuremath{\bm{f}}}

\newcommand{\bk}{\ensuremath{\bm{k}}}

\newcommand{\bp}{\ensuremath{\bm{p}}}
\newcommand{\bq}{\ensuremath{\bm{q}}}
\newcommand{\br}{\ensuremath{\bm{r}}}

\newcommand{\bu}{\ensuremath{\bm{u}}}

\newcommand{\bx}{\ensuremath{\bm{x}}}
\newcommand{\by}{\ensuremath{\bm{y}}}

\newcommand{\bphi}{\ensuremath{\bm{\phi}}}
\newcommand{\bone}{\ensuremath{\bm{1}}}

\newcommand{\Fcal}{\ensuremath{\mathcal{F}}}
\newcommand{\Ncal}{\ensuremath{\mathcal{N}}}
\newcommand{\Tcal}{\ensuremath{{\mathcal{T}}}}
\newcommand{\Ocal}{\ensuremath{{\mathcal{O}}}}
\newcommand{\Mcal}{\ensuremath{{\mathcal{M}}}}
\newcommand{\Acal}{\ensuremath{{\mathcal{A}}}}
\newcommand{\Kcal}{\ensuremath{{\mathcal{K}}}}
\newcommand{\qbb}{\ensuremath{{\mathcal{Q}}}}
\newcommand{\phibb}{\ensuremath{{\mathcal{T}}}}
\newcommand{\Rcal}{\ensuremath{\mathcal{R}}}

\newcommand{\MLOP}{\ensuremath{\Tcal_{M2L}}}

\definecolor{orange}{cmyk}{0,0.353,1,0}

\theoremstyle{definition}

\newtheorem{alg}{Algorithm}

\journal{Journal of Computational Physics}

\begin{document}

\begin{frontmatter}


 \author[add1,add2]{Wen Yan\corref{cor1} }
 \ead{wyan@flatironinstitute.org}
 \ead{wenyan4work@gmail.com}
 \author[add1,add2]{Michael Shelley }
 \ead{mshelley@flatironinstitute.org}
\cortext[cor1]{Corresponding author}


\title{Flexibly imposing periodicity in kernel independent FMM: A Multipole-To-Local operator approach}

\address[add1]{Center for Computational Biology, Flatiron Institute, Simons Foundation, New York 10010}
\address[add2]{Courant Institute of Mathematical Sciences, New York University, New York 10010}

\begin{abstract}
An important but missing component in the application of the kernel independent fast multipole method (KIFMM) is the capability for flexibly and efficiently imposing singly, doubly, and triply periodic boundary conditions. In most popular packages such periodicities are imposed with the hierarchical repetition of periodic boxes, which may give an incorrect answer due to the conditional convergence of some kernel sums. Here we present an efficient method to properly impose periodic boundary conditions using a near-far splitting scheme. The near-field contribution is directly calculated with the KIFMM method, while the far-field contribution is calculated with a multipole-to-local (M2L) operator which is independent of the source and target point distribution. The M2L operator is constructed with the far-field portion of the kernel function to generate the far-field contribution with the downward equivalent source points in KIFMM. This method guarantees the sum of the near-field \& far-field converge pointwise to results satisfying periodicity and compatibility conditions. The computational cost of the far-field calculation observes the same $\Ocal(N)$ complexity as FMM and is designed to be small by reusing the data computed by KIFMM for the near-field. The far-field calculations require no additional control parameters, and observes the same theoretical error bound as KIFMM. We present accuracy and timing test results for the Laplace kernel in singly periodic domains and the Stokes velocity kernel in doubly and triply periodic domains.
\end{abstract}


\end{frontmatter}


\section{Introduction}
\label{sec:intro}
Since its invention in the 1980s \cite{greengard_fast_1987}, the Fast
Multipole Method (FMM) has been applied in many areas of computational
science, with its popularity due to its $\Ocal(N)$ complexity, efficient
arithmetic intensity in implementation \cite{yokota_fmm_2012}, and
better scalability compared to competing methods like FFT in perfoming
the N-body sum:
\begin{align}
	q_i^t = \sum_{j} K(\bx_i,\by_j) \phi_j^s, \quad \forall i,j\in \{1,2,...,N\},
\end{align}
where $q_i^t$ is the potential located at $\bx_i$, $\phi_j^s$ is the source intensity  located at $\by_j$, and $K(\bx_i,\by_j)$ is  the kernel function. The original version of FMM utilizes multipole expansions, and triple
periodicity (TP) can be implemented with an explicit Ewald sum of the multipole basis functions \cite{kudin_revisiting_2004}. However, performing Ewald sums on-the-fly is expensive and some FMM packages simply approximate periodic boundary conditions by including the nearest
image of the original box \cite{yokota_petascale_2013}. Also, single and double periodicity (SP and DP) implementations are not as
straightforward as the TP case because the Ewald summation of the basis functions may involve special functions which are difficult to
evaluate to high precision  \cite{tornberg_ewald_2015,harris_incomplete_2008}.

Kabadshow \cite{kabadshow_periodic_2012} proposed
a periodizing method based on splitting the near and far-field
contributions of the Laplace kernel. Contributions from periodic
images adjacent to the original box $B_0$ are considered `near-field',
and contributions from all other periodic image boxes are considered
`far-field'. The far-field contributions are from boxes well-separated
from $B_0$ and can be approximated by multipole expansion. The
far-field contributions are precomputed to high precision as a
multipole-to-local (M2L) operator applied on multipole basis
functions.  With this scheme, SP, DP and TP boundary conditions can be
calculated efficiently without invoking Ewald after the precomputing
stage, as long as the corresponding M2L operators are properly
evaluated. However, this is an explicit \emph{summation} method which
relies on a specific order of summing the periodic images of the
multipole basis functions, and is hard to generalize to kernel
independent FMM in the absence of multipole basis functions.

The kernel independent fast multipole method
(KIFMM) \cite{ying_kernel-independent_2004} reformulated the FMM into a
`method of fundamental solutions' (MFS) form, by replacing the multipole
expansion with a set of `equivalent sources'. In this way, the same
code can be used for different kernel functions $K(\bx,\by)$. 
Some recent implementations of KIFMM have used a Hierarchical Repetition (HR)
method, to hierarchically `copy-and-paste'
the original box into a rectangular tiling with $(N_x,N_y,N_z)$ periodic images in $(x,y,z)$ directions, where $(N_x,N_y,N_z)$ is usually taken to be a cube. The copy
is hierarchical in the sense that, for a $k$ level hierarchy,
$N_x\sim2^k$. This method reuses the adaptive octree (or quadtree in
2D) structure built for the original box $B_0$, and is very efficient
because the hierarchical repetition sum of all equivalent sources can
be precomputed and stored for reuse as an M2L translation operator.
However, this method does not
always give the correct answer for typical kernels with $~1/r$ decay,
like the Laplace or the Stokes kernel. This is because the periodic
sum of such $1/r$ kernels is conditionally convergent, which makes their
results subject to the order of summation.For example, a cubic sum
with $N_x=N_y=N_z$ and a rectangular sum $N_x=2N_y=2N_z$ give
different answers due to the conditional convergence, and both results are not correct. A recent implementation \cite{malhotra_pvfmm_2015} of HR method resolved this conditional convergence by removing the net monopole and net dipole moments in the original box $B_0$ in the hierarchical sum. With this method, finite results can be achieved but the physical compatibility conditions are often not satisfied. In most physical applications the conditional convergence
problem is resolved with the well-known Ewald
sum \cite{tornberg_ewald_2015,lindbo_fast_2011}. The Ewald sum imposes
the exact periodicity by starting from a Fourier expansion of the
field in the original box $B_0$, and therefore is free from the
conditional convergence error arising from particular orders of
summation.

Recently, methods have been proposed to circumvent the conditional
convergence problem for both FMM and KIFMM by directly imposing
periodicity on a set of check points placed on a chosen
surface \cite{barnett_unified_2016,gumerov_method_2014}.  These methods
also follow the idea of near-far splitting, where the near-field is
directly calculated by FMM, and the far-field is \emph{solved for} on a set of `equivalent sources' instead of being
\emph{summed} as developed by
Kabadshow \cite{kabadshow_periodic_2012}. The
locations of the equivalent sources are usually on a surface outside
the original box, and are chosen according to a high-accuracy
quadrature rule suitable for the dimension and geometry. Gumerov and
Duraiswami \cite{gumerov_method_2014} demonstrated the application of this idea in the triply periodic
electrostatic problem in 3D space. In more recent work, Barnett et
al. \cite{barnett_unified_2016} systematically analyzed the method in a
doubly periodic domain for both the Laplace and Stokes
kernels in 2D space. However, partial periodicities (SP and DP in 3D space) are
not easy to impose because imposing the partial periodicity on the
check points is sometimes not sufficient to guarantee the periodicity
in the entire periodic box, and the `zero in infinity' condition in non-periodic directions must be supplied to determine the solution.

In this work, we develop a simple method to compute  the periodic
kernel independent FMM for singly, doubly and triply periodic boundary
conditions. The method follows the idea of near-far splitting and the M2L
operator, but keeps the KIFMM formulation. The method is an
improvement of the methods discussed above and is designed to:
\begin{itemize}
	\item Reuse the data already calculated by KIFMM.
	\item Minimize the modification of the underlying KIFMM code.
	\item Converge close to machine precision with known error analysis.
	\item Keep the $\Ocal(N)$ complexity of KIFMM and minimize \verb|MPI| communication.
	\item Add no extra control points, equivalent sources, or
          tweaking parameters to KIFMM.
\end{itemize}
We describe our general method in Section~\ref{sec:method}, and
present accuracy and timing results for the Laplace kernel in SP geometry
and the Stokes velocity kernel in DP and TP geometries in
Section~\ref{sec:results}. Finally in Section~\ref{sec:discuss} we
discuss the possible optimization, extensions and applications of this
method. 


\section{Methodology}
\label{sec:method}
\subsection{The general idea}

Before diving into the formulation of KIFMM, we first describe the general idea of our algorithm by comparing it with previous work \cite{barnett_unified_2016,gumerov_method_2014}, in which the periodicity is \emph{solved} from a linear system. As a simple example, consider a unit rectangle $[0,1)^2$ in 2D space, with its left, right, top and bottom boundary denoted as $L,R,T,B$. A doubly periodic function $u(x,y)$ satisfies a Poisson equation with sources $\phi_k$ located at $\br_k$ and periodic boundary conditions:
\begin{align}\label{eq:poisson2d}
	\nabla^2 u &=\sum \phi_k\delta(\br-\br_k), \\
	u(x+1,y) &= u(x,y),\quad \forall x,y, \\
	u(x,y+1) &= u(x,y),\quad \forall x,y.
\end{align}
Barnett et al. \cite{barnett_unified_2016} showed that the solution to this equation is unique up to a constant. This constant can be later determined by some other physical information of the problem, like a fixed temperature at some point. Therefore after adding the equivalent sources to the system, the strength of the equivalent sources can be determined uniquely to approximate the far field. Gumerov and Duraiswami used a similar formulation in 3D, except that the periodicity is checked for points on a spherical (in 3D) surface instead of on the $L,R,T,B$ boundaries of the domain. Similar uniqueness argument is also proved in the work by Barnett et al. \cite{barnett_unified_2016} for doubly periodic Stokes problems on 2D domains.

When the above method is extended to partial periodicity, uniqueness is no longer guaranteed by simply imposing the periodic boundary condition. When periodicity in the $y$ direction is removed from the above equation, we have a general solution to the 2D Laplace equation with single periodicity in $x$-direction with any integer $k$:
\begin{align}\label{eq:partialperiodicgensol}
	u^{SP} &= a_0 + b_0 y + \sum_{k\neq 0} c_k \exp \left(2k\pi i x\right)\exp\left(2k\pi y\right). 
\end{align}
With any real $a_0,b_0$ and $c_k$, $u^{SP}$ satisfies the partially periodic Laplace equation:
\begin{align}
	\nabla^2 u =0,\quad u(x+1,y) &= u(x,y),\quad \forall x,y.
\end{align}
For any solution $u$ to Eq.~(\ref{eq:poisson2d}) with periodicity in $x$ direction only, $u+u^{SP}$ is still a solution. 
Therefore, extra conditions must be supplied in the $y$ direction to solve for the constants $a_0,b_0$ and $c_k$ to maintain the uniqueness of the above method, otherwise the strength of the added equivalent sources cannot be determined. For Helmholtz kernels such decaying conditions can be determined as $y\to\pm\infty$ as demonstrated in \cite{cho_robust_2015}, and systems with Laplace/Stokes kernels could be similarly solved. Similar families of $u^{SP}$ are also straightforward to construct for the Stokes equations by taking a Fourier expansion in the periodic directions Eq.~(\ref{eq:partialperiodicgensol}).

Instead of imposing periodicity by \emph{solving} the governing PDEs with periodicities, which relies on the uniqueness of solutions, we choose to directly start from the periodic Green's function $K^P$ for those PDEs. Those periodic Green's functions are usually explicitly known as absolutely convergent series, because they are usually derived by the Fourier analysis of the corresponding PDE. With those $K^P$, the following summation automatically satisfy the conditions in both periodic and non-periodic directions:
\begin{align}\label{eq:KPsum}
q_i^t = \sum_{j} K^P(\bx_i,\by_j) \phi_j^s, \quad \forall i,j\in \{1,2,...,N\}.
\end{align}
However evaluating $K^P$ for every pair of source-target is usually expensive, and our objective is to use the freespace kernel $K$ with KIFMM, and to ensure the results converge pointwise to Eq.~(\ref{eq:KPsum}). To achieve this, we split the periodic infinite domain into a small near field and a far field. The near field is easy to calculate with KIFMM because it is simply a finite system with open boundary conditions. The effect of periodicity is added by letting the far field generate the remaining piece of a true periodic solution. The contributions from far field are smooth on the original domain, which allows us to approximate the solution with a few equivalent sources, following the idea of KIFMM.

To apply the method presented in this work to a new kernel $K$, we need to know its periodic form $K^P$ explicitly. In this work $K^P$ is constructed with Ewald methods for convenience, because Ewald methods have been extensively used for simulations of Laplace and Stokes systems. However, the Ewald methods are not the only choice. $K^P$ can also be constructed by direct sums in either real space or Fourier space if the convergence is rapid. For example for the fast decaying Yukawa potential $\exp(-\kappa r)/r$ with a real $\kappa$, a direct summation of several image boxes in real space is often sufficient (depending on the value of $\kappa$). For the biharmonic equation in 3D space, a series summation in Fourier space is also sometimes sufficient because the Green's function converges as $1/k^4$ in Fourier space. If $K^P$ is not easy to find analytically, the method of Barnett et al. \cite{barnett_unified_2016} can be used if the uniqueness is guaranteed as discussed above. The Hierarchical Repetition method in \cite{malhotra_pvfmm_2015} is also applicable if the far-field is properly fixed to satisfy the physical compatibility condition. The central idea of our method is to construct the near \& far field converging pointwise to Eq.~(\ref{eq:KPsum}) when summed together so that the periodicity and physical compatibility condition are automatically satisfied by $K^P$. We do so efficiently by reusing the KIFMM data for near field to construct the far field.

\subsection{Formulation and algorithm}
Consider the following formulation of the periodic FMM problem:

\textbf{Problem.} Given $T$ target points, $S$ source points each with strength
        $\phi^s$ in a box $B_0=[0,L_x)\times[0,L_y)\times[0,L_z) $,
              and a kernel function $K(\bx_{t},\by_{s})$, evaluate the
              potential on each target point with periodic BC:
	\begin{equation}\label{eq:pfmmDef}
		q^t = \sum_{\bp\in \bP} K(\bx_{t},\by_{s}+\bp) \phi^s.
	\end{equation}
Here $t\leq T$ and $ s\leq S$ are the indices for target and source
points, $\bx_t$ is the location of target points, $\by_s$ is the
location of source points, and repeated the index $s$ implies a
summation over all
$s\in S$. $\bP$ is the set of all periodic vectors. For example, for a
unit cubic box periodic in $z$ direction: $\bP = \{\bp\in \mathbb{Z}^3
: p_x=0,p_y=0,p_z\in \mathbb{Z}\}$.

We will assume for simplicity that all target and source points are
distributed in a unit cubic box $B_0=[0,1)^3$, which we will call the
  ``original box''. All other periodic boxes (with $\bp\neq0$) are
  called ``image boxes''.

We follow the original KIFMM notations \cite{ying_kernel-independent_2004}:

\begin{tabular}{ll}
$B_i$ 	& a box in the $i$th level of the octree \\
$\Ncal^B$	& the near-field of $B$, including $B$ itself \\
$\Fcal^B=\Rcal^d/\Ncal^B$ & the far-field of $B$, where $d$ is the spatial dimension\\
$I_s^B$ & the set of indices of source points inside $B$\\
$I_t^B$ & the set of indices of target points inside $B$\\
$\phi^{B,s}$ & the strength of source points inside $B$\\
$q^{B,t}$ & the potential at target points inside $B$\\
$\by^{B,u}$ & the upward equivalent surface of $B$\\
$\phi^{B,u}$ & the upward equivalent density of $B$\\
$\bx^{b,u}$ & the upward check surface of $B$\\
$q^{B,u}$ & the upward check potential of $B$\\
$\by^{B,d}$ & the downward equivalent surface of $B$\\
$\phi^{B,d}$ & the downward equivalent density of $B$\\
$\bx^{B,d}$ & the downward check surface of $B$\\
$q^{B,d}$ & the downward check potential of $B$\\
$p$ & the number of grid points per cube edge for the discretization of equivalent surfaces\\
$s$ & the maximum number of source (or target) points allowed in a leaf box\\
$N$ & the total number of source and target points\\
\end{tabular}\newline

\begin{figure}

\centering \includegraphics[width=0.75\linewidth]{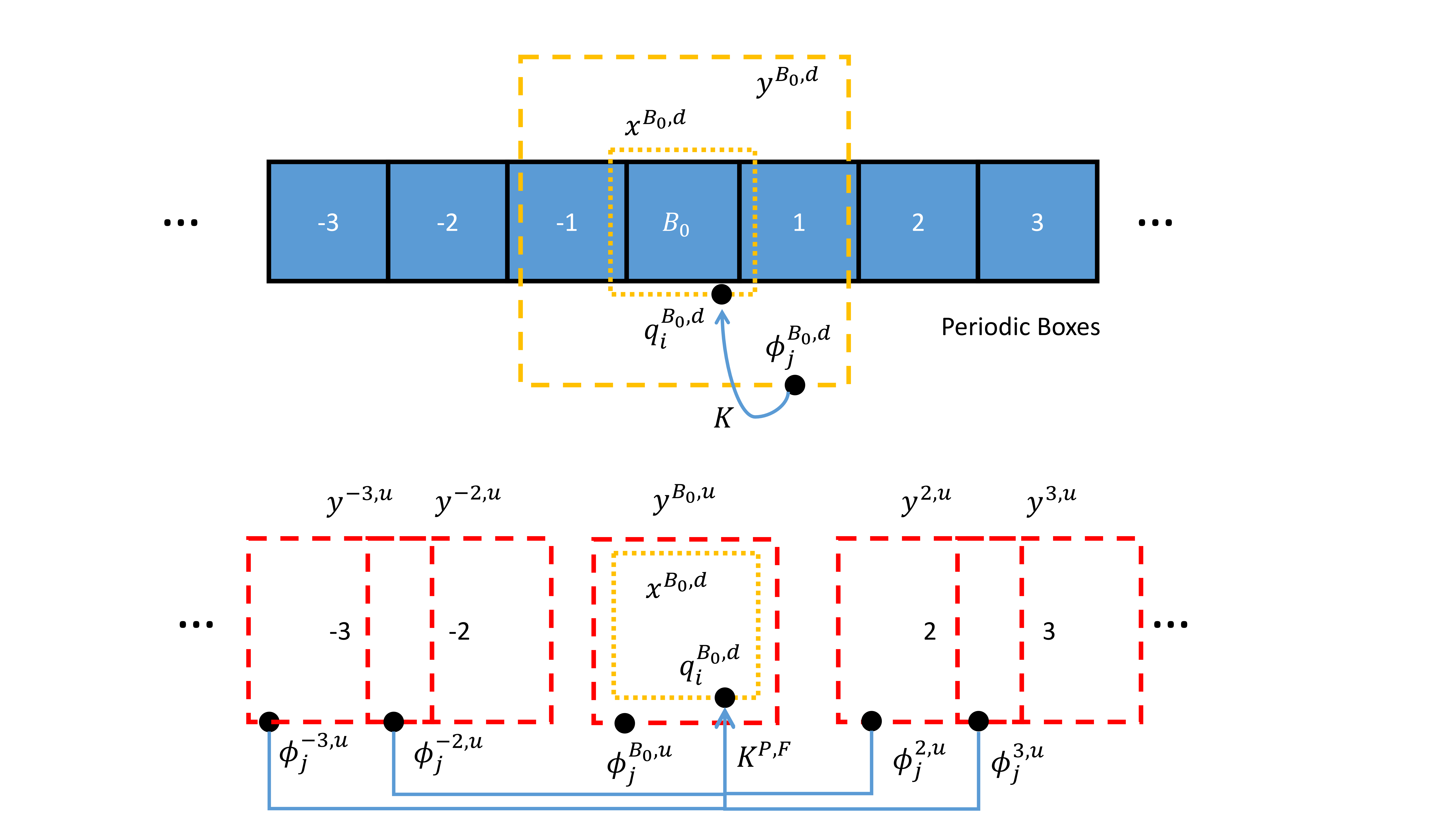}
\caption{A graphical representation of Algorithm~\ref{alg2} for the singly
  periodic boundary condition. $\Ncal^{B_0}$ is chosen with $\ell=1$
  layer of image boxes. The top figure shows the downward check and
  equivalent surfaces, and the bottom figure shows the downward check
  surface of $B_0$ and all upward equivalent surfaces of periodic
  image boxes belonging to $\Fcal^{B_0}$. The algorithm is to find a set
  of (discretized) downward equivalent sources $\phi_j^{B_0,d}$ to
  match the field in $B_0$ generated by all
  $\phi_j^{u}\in\Fcal^{B_0}$. Imposing periodicity gives all
  $\phi_j^{u}\in\Fcal^{B_0}$ equal to $\phi_j^{B_0,u}$, which is
  calculated by FMM for $\Ncal^{B_0}$. In this figure the surfaces $\by^{B_0,u}$ and $\bx^{B_0,d}$ are plotted as not overlapping with each other, but in real calculations they are usually chosen to be identical. }
\label{fig:schematic}
\end{figure}

We choose to follow the original KIFMM
notation \cite{ying_kernel-independent_2004} to ease the
comparison. This is in contrast to more standard notations, say where
$q$ means charge and $\phi$ means potential. For each cubic box in the octree, its equivalent and check surfaces are cubic surfaces enclosing the octree box as defined in \cite{ying_kernel-independent_2004}. The equivalent and check surfaces for the root box are illustrated in Figure~\ref{fig:schematic}. The discretization in 3D
space is a regular grid on the cube surface of equivalent densities,
and $p$ includes the cube vertices. $\Ncal^B$ includes boxes adjacent
to $B$, within a distance of $\ell$ times $B$'s edge length. Following the definition of Ying et al \cite{ying_kernel-independent_2004}, $\Ncal^B$ includes boxes in the same level of $B$ in the octree. For example, in this work $B_0=[0,1)^3$ refers to the unit cubic original
box, and $\Ncal^{B_0}(\ell=1)$ includes 2 neighboring image boxes in
  a singly periodic geometry and 26 neighboring image boxes in a
  triply periodic geometry. In the case of $\ell=2$,
  $\Ncal^{B_0}(\ell=2)$ includes 4 and 124 neighboring boxes in those
  two cases respectively.  Depending on the boundary conditions,  the
  number of neighboring boxes in $\Ncal^{B}$ may be different from
  $\Ncal^{B_0}$ for a leaf box $B$ in the octree. For example, in a doubly periodic system $\Ncal^{B}$ is smaller than $\Ncal^{B_0}$, for a leaf box $B$ at the open boundary. 
  
The fundamental idea in KIFMM is that for a set of target points in
$B$, the field due to source points far away (from $\Fcal^B$) is
smooth enough to be approximated by an equivalent source surface
enclosing the target points. The strength distribution on the
equivalent surface can be calculated by matching the field strength
from the source points and from the equivalent surface on another
check surface. The density $\phi^{B,u}$ on $\by^{B,u}$
approximates the effect of source points in the box $B$ ($\phi^{B,s}$
with $s\in I_s^B$) to all far away target points, and is checked on
$\bx^{b,u}$. The density $\phi^{B,d}$ on $\by^{B,d}$
approximates the effect of all far away source points on the potential
of target points in the box $B$ ($q^{B,t}$ with $t\in I_t^B$ ), and is
checked on $\bx^{b,d}$. Application of this formulation
within an adaptive octree structure leads to the following five
crucial steps of KIFMM:

\begin{itemize}
	\item \makebox[3em]{S2M:\hfill} The \emph{Source to Multipole} operation evaluates $\phi^{B,u}$ with $\phi^{B,s}$ in a leaf box $B$.
	\item \makebox[3em]{M2M:\hfill} The \emph{Multipole to Multipole} operation transforms $\phi^{B,u}$ of a box's children to its own $\phi^{B,u}$.
	\item \makebox[3em]{M2L:\hfill} The \emph{Multipole to Local} operation transforms $\phi^{B,u}$ of a box to the $\phi^{B,d}$ of a non-adjacent box.
	\item \makebox[3em]{L2L:\hfill} The \emph{Local to Local} operation transforms the $\phi^{B,d}$ of a box's parent to its own $\phi^{B,d}$.
	\item \makebox[3em]{L2T:\hfill} The \emph{Local to Target}
          operation evaluates $q^{B,t}$ with known $\phi^{B,d}$.
\end{itemize} 
Besides these five steps, for leaf boxes not well separated, S2T operations are directly applied on every pair of source and target points in those boxes to calculate the contributions to $q^t$.

We follow the idea of KIFMM, and split the periodic geometry into the near-field $\Ncal^{B_0}$ and the far-field $\Fcal^{B_0}$, with an adjustable splitting layer number $\ell$:
\begin{align}\label{eq:pfmmSplit}
	q^t(\bx^t\in B_0) &= \underbrace{\sum_{ \bp\in \Ncal^{B_0}} K(\bx_{t},\by_{s}+\bp)\phi^s }_{\text{near-field: } q_{\Ncal}^t} + \underbrace{\sum_{ \bp\in \Fcal^{B_0}} K(\bx_{t},\by_{s}+\bp)\phi^s}_{\text{far-field: } q_{\Fcal}^t},
\end{align}
where $q_{\Ncal}^t$ is a finite sum for a well defined `free-space'
kernel function $K$, and is straightforward to calculate with
KIFMM. As demonstrated in Section~\ref{subsec:effectl} we found that adjusting $\ell$ has almost no effect on the
accuracy and timing of the algorithm, and so we fix $\ell=2$ for
results in this work. Also, due to the periodic structure,
the octree and multipoles can be built and calculated for $B_0$ only,
and then we can `copy-and-paste' to include all image boxes in
$\Ncal^{B_0}$. This costs much less than directly building an octree and
evaluating multipoles for the extended domain $\Ncal^{B_0}$, and this
`copy-and-paste' method is usually already implemented in most KIFMM
packages performing hierarchical repetition calculations
\footnote{Some minor modifications to the hierarchical repetition
  routines may be necessary because $\Ncal^{B_0}$ always includes an
  integer $\ell$ of images of $B_0$ in every periodic direction,
  instead of possibly fractional copies in the hierarchical repetition
  routines.}.

Because $\Fcal^{B_0}$ is well separated from $B_0$, the far-field
contribution $q_{\Fcal}^t$ can be approximated by the equivalent
densities\footnote{In the following $K(\bx_t,\by^u)\phi^u$ and
  $K(\bx_t,\by^d)\phi^d$ should be understood as an integration over
  the surfaces $\by^{B,u}$ and $\by^{B,d}$.}:
\begin{align}\label{eq:qtf1}
	q_{\Fcal}^t = \sum_{\bp\in \Fcal^{B_0}} K(\bx_{t},\by^{B_0+\bp,u})\phi^{B_0+\bp,u} = \sum_{\bp\in \Fcal^{B_0}} K(\bx_{t},\by^{B_0+\bp,u})\phi^{B_0,u},
\end{align}
where we utilized the periodicity of $B_0$:
$\phi^{B_0+\bp,u}=\phi^{B_0,u}$, and $\phi^{B_0,u}$ is already
calculated in the KIFMM calculations for $q_{\Ncal}^t$. The far-field
sum can then be written as an `exact periodic' part and an
`near-field' part:
\begin{align}\label{eq:qtf2}
q_{\Fcal}^t &= \sum_{\bp\in \bP} K(\bx_{t},\by^{B_0+\bp,u})\phi^{B_0,u} - \sum_{\bp\in \Ncal^{B_0}} K(\bx_{t},\by^{B_0+\bp,u})\phi^{B_0,u}, \\
	&=K^P(\bx_{t},\by^{B_0,u})\phi^{B_0,u} - \sum_{\bp\in \Ncal^{B_0}} K(\bx_{t},\by^{B_0+\bp,u})\phi^{B_0,u},
\end{align} 
where the periodic kernel $K^P$ is defined with the free-space kernel $K$:
\begin{align}
	K^P = \sum_{\bp\in \bP } K(\bx_{t},\by_{s}+\bp).
\end{align} 
We can further define a far-field periodic kernel:
\begin{align}
K^{P,F}(\bx_{t},\by_{s}) = K^P(\bx_{t},\by_{s}) - \sum_{\bp\in \Ncal^{B_0}}  K(\bx_{t},\by_{s}+\bp),
\end{align}
and the far-field contribution is simplified:
\begin{align}
q_{\Fcal}^t = K^{P,F}(\bx_{t},\by^{B_0,u})\phi^{B_0,u}.
\end{align}

For the Laplace kernel $K=1/\ABS{\bx_{t}-\by_{s}}$, $K^P$ is
well-known as the Ewald summation, and its singly, doubly and triply
periodic formulations in 3D space are given by
Tornberg \cite{tornberg_ewald_2015}. For Stokes kernels for velocity,
pressure, stress, etc., the Ewald sums for $K^P$ are also well-known
 \cite{lindbo_fast_2011,lindbo_spectrally_2010,sierou_accelerated_2001,sierou_accelerated_2002,wang_spectral_2016}
and can also be calculated by a transform from the Laplace kernel to
the Stokes kernel \cite{tornberg_fast_2008}. For more general form of
kernels, they can usually be expressed by a combination of $1/r^n$
with different $n$, or as a sum of spherical multipole basis
functions \cite{schmidt_multipole_1997,mazars_ewald_2010,yan_behavior_2016}.

The existence of $K^P$ usually depends on some compatibility
conditions, related to the physical setting. Such a compatibility
condition usually manifests itself as requiring the convergence
of $K^P$. For example, for an electrostatics problem the net
charge in the original box $B_0$ must be zero, otherwise the periodic
sum of potentials is divergent. For Stokes problems, the compatibility
condition is different for different periodic geometries, and will be
discussed later in Section~\ref{sec:results}. In general, the
compatibility condition usually takes either (or both) of the two following
forms:
\begin{align}\label{eq:compatcondition}
  &\text{Neutrality: }&\sum_{s\in S} \phi^s &=0,\\
  &\text{Specified net field integration: }& \int_{B_0}
  q^t(\bx) d^3\bx&=C,
\end{align}
where $C$ is a given constant and is usually $0$. In some cases the above two conditions may be modified to fit the physical setting, and one example is given in Section~\ref{subsec:stokes2P}.

With known $K^P$ for a given $K$, we can evaluate the periodic FMM
with the following algorithm:
\begin{alg}\label{alg1}
	Evaluate periodic FMM~(\ref{eq:pfmmDef}) with near-far splitting method
	\begin{enumerate}[a).]
		\item Call KIFMM to calculate $q_{\Ncal}^t$.
		\item Read $\phi^{B_0,u}$ from the KIFMM routine.
		\item Calculate $q_{\Fcal}^t$ with Eq.~(\ref{eq:qtf2}) and add it to $q_{\Ncal}^t$.
	\end{enumerate}
\end{alg} 
This algorithm is straightforward to implement, but is far from
optimal because the evaluation of $q_{\Fcal}^t$ requires frequent
calls to $K^P$, which is usually very expensive although in practice
$K^P$ is always calculated with some accelerated methods like Particle
Mesh Ewald or Spectral
Ewald \cite{lindbo_spectrally_2010,wang_spectral_2016}.

In fact, step c) can be accelerated with the M2L translation
operation, as KIFMM does for all well separated boxes. The operation
is to solve a first-kind integral equation on the downward equivalent
and check surfaces, to find the
equivalent density generating fields in $B_0$ matching the periodic
image boxes in $\Fcal^{B_0}$:
\begin{align}\label{eq:M2Lsolve}
\forall \bx\in\bx^{B_0,d}:	\int_{\by^{B_0,d}} K(\bx,\by) \phi^{B_0,d}(\by) d\by = \int_{\by^{B_0,u}}K^{P,F}(\bx,\by^{B_0,u})\phi^{B_0,u} d\by.
\end{align}

With the discretization of equivalent and check surfaces, we have a linear equation:
\begin{align}
	\Acal \phi^{B_0,d} = \Kcal^{P,F} \phi^{B_0,u},
\end{align}
where the matrices $\Acal$ and $\Kcal^{P,F}$ will be explicitly demonstrated in the next section.
The solution is a \emph{linear} M2L translation operator $\MLOP=\Acal^\dagger \Kcal^{P,F}$ on arbitrary $\phi^{B_0,u}$:
\begin{align}
	\phi^{B_0,d}(\by) = \MLOP \phi^{B_0,u}.
\end{align}
Here we use the pseudo-inverse $\Acal^\dagger$ because $\Acal$ would be numerically usually singular. We discuss this further in Section~\ref{subsec:SVD}.

This operator depends only on $K$, $K^{P}$, and the check and
equivalent surfaces, and is independent of the source and target point
distributions in the box $B_0$. With given locations and
discretizations of the upward and downward equivalent surfaces,
$\MLOP$ need be calculated only once for all simulations. With this method,
we propose an algorithm much faster than Algorithm~\ref{alg1}:
\begin{alg}\label{alg2}
	Evaluate periodic FMM~(\ref{eq:pfmmDef}) with a near-far
        splitting method and $\MLOP$ operator\\
        Stage 1. Precomputing.
	\begin{enumerate}[a).]
		\item Calculate $\MLOP$ by solving Eq.~(\ref{eq:M2Lsolve}).
	\end{enumerate}
        Stage 2. KIFMM Evaluation.
	\begin{enumerate}[a).]
		\item Call KIFMM to calculate $q_{\Ncal}^t$.
		\item Read $\phi^{B_0,u}$ from the KIFMM routine.
		\item Calculate $\phi^{B_0,d}(\by) = \MLOP \phi^{B_0,u}$.
		\item Evaluate $q_{\Fcal}^t=K(\bx_{t},\by^{B_0,d})\phi^{B_0,d}$
                  and add it to $q_{\Ncal}^t$.
	\end{enumerate}
\end{alg}
Stage 1 is illustrated in Figure~\ref{fig:schematic}. This is much
faster than Algorithm~\ref{alg1} because it involves only the
free-space kernel $K$ in the evaluation stage.
The two stages of Algorithm~\ref{alg2} are similar to the method
proposed by Gumerov and Duraiswami \cite{gumerov_method_2014}, but
Algorithm~\ref{alg2} is much faster because it requires no extra
control points and reuses the octree data for near-field
calculations. In fact Algorithm~\ref{alg2} can be further accelrated with proper modifications of the underlying KIFMM code, by replacing step d) with a downward pass of L2L and L2T operations through the octree. However as shown in Section~\ref{subsec:stokestiming} the cost of step d) is already small and we choose to keep the underlying KIFMM code as simple as possible, without applying this optimization.

The discretization and locations of the check and equivalent surfaces
are the key component of this method. In this work we follow the
choice of the original KIFMM method \cite{ying_kernel-independent_2004}
and the high-performance KIFMM package
PVFMM \cite{malhotra_pvfmm_2015,malhotra_algorithm_2016}. For the
original box $B_0=[0,1)^3$, the locations of $\by^{B,u}$ and
  $\bx^{B,d}$ are identical, with both being a cube with edge length
  $1.05$ centered at $(0.5,0.5,0.5)$. The locations of $\bx^{B,u}$ and
  $\by^{B,d}$ are also identical, being a cube with edge length $2.95$
  also centered at $(0.5,0.5,0.5)$. The discretization is chosen to be
  a regular mesh uniformly distributed on the 6 surfaces of a cube. On
  each edge $p$ points (including two vertices) are uniformly
  distributed, and in total $6(p-1)^2+2$ points are distributed on the
  cube surface.

\subsection{The method to solve for $\MLOP$}
The theory about approximating $\MLOP$ with Eq.~(\ref{eq:M2Lsolve}), with
rigorous error analysis, is well understood in terms of both integral
equation theory and for the KIFMM method
 \cite{ying_kernel-independent_2004}. In the following we briefly
describe the structure of $\MLOP$ and the method to solve for it. For a general kernel function $\bK: \Rcal^{k_s}\rightarrow\Rcal^{k_t} $, $k_s$ is usually termed as the source dimension of $\bK$, and $k_t$ is usually termed as the target dimension. For
example, the Laplace kernel is defined on $\Rcal\rightarrow\Rcal$, the Stokes velocity kernel on $\Rcal^3\rightarrow\Rcal^3$, and the Stokes pressure kernel on
$\Rcal^3\rightarrow\Rcal$. With a
chosen discretization of $p$ points per cube edge,
the operator $\MLOP$ takes a block column form:
\begin{align}
\MLOP = \left[\phibb_{1}^{B_0,d},\phibb_{2}^{B_0,d},\cdots,\phibb_{n}^{B_0,d}\right]\in \Rcal^{k_s n \times k_s n},
\end{align}
where $n=6(p-1)^2+2$ is the total number of discretization points on the equivalent surface. Each $\phibb_{j}^{B_0,d},1\leq j \leq n$, has a block row structure:
\begin{align}
\phibb_{j}^{B_0,d}=\left[ \bphi_{1j}^{B_0,d}, \bphi_{2j}^{B_0,d}, \cdots,  \bphi_{nj}^{B_0,d} \right]^T ,
\end{align}
where each small square block $\bphi_{ij}^{B_0,d}\in\Rcal^{k_s\times k_s}$, and denotes the strength of a source needed at location
$i$ on the downward equivalent surface of $B_0$, due to the source with unit strength at
location $j$ on the upward equivalent surface of $B_0$.

Each block column vector $\phibb_{j}^{B_0,d}$ is solved for through a
linear equation,
\begin{align}\label{eq:MLOPAxb}
\Acal \phibb_{j}^{B_0,d} = \qbb_j^{B_0,d},
\end{align}
where $\qbb_j^{B_0,d}$ also has a row block structure:
\begin{align}
\qbb_{j}^{B_0,d}&=\left[ \bq_{1j}^{B_0,d}, \bq_{2j}^{B_0,d}, \cdots,  \bq_{nj}^{B_0,d}  \right]^T,\\
\forall i,\quad \bq_{ij}^{B_0,d}& = \bK^{P,F} \left(\bx_{i}^{B_0,d},\by_{j}^{B_0,u}\right).
\end{align}
Each block $\bq_{ij}^{B_0,d}$, having dimension $(k_t,k_s)$, denotes
the far-field potential induced on check point $i$
by equivalent source $j$. It is worth noting that while the
upward equivalent surface $\by^{B_0,u}$ is chosen to overlap with the
downward check surface $\bx^{B_0,d}$, the blocks $\bq_{ii}^{B_0,d}$ (here the repeated index $i$ does not mean summation) are not singular because $\bK^{P,F}$ sums only the far-field contributions from $\Fcal^{B_0}$.

The matrix $\Acal$ in Eq.~(\ref{eq:MLOPAxb}) also has a block structure based on the free-space kernel $\bK$:
\begin{align}
\Acal_{ij} = \bK\left( \bx_{i}^{B_0,d},\by_{j}^{B_0,d} \right).
\end{align}
Usually the check and equivalent surfaces, $\bx^{B_0,d}$ and
$\by^{B_0,d}$, are discretized with the same discretization parameter $p$
in KIFMM implementations, and $\Acal$ is block square. Some work in
the literature suggests using a finer discretization on
the check surface and performing least
squares \cite{gumerov_method_2014} to solve for $\qbb_j^{B_0,d}$ instead of solving Eq.~(\ref{eq:MLOPAxb}). We have found no benefit in accuracy from oversampling, and always use the same $p$ for check and equivalent surfaces, at least for the
results reported in this work. For kernels with dimension $k_t<k_s$,
oversampling may be necessary to maintain the square shape of
$\Acal$. It is well-known that $\Acal$ can be numerically nearly singular with large condition number \cite{barnett_stability_2008}, and
we shall discuss briefly the method to solve for it in
Section~\ref{subsec:SVD}.

\subsection{Cost and error}
With a chosen discretization parameter $p$, the cost of each stage of
Algorithm~\ref{alg2} can be estimated. In the precomputing stage,
$\MLOP$ is a matrix of $\Ocal((6(p-1)^2+2)^2)\sim\Ocal(p^4)$ entries,
and solving for it requires the pseudo-inverse of a matrix of the same
dimension. The theoretical cost of forming pseudo-inverses is
well-known and we skip its discussion. The computation usually completes within $100$ seconds and uses less than $1GB$ of memory, and we give more details about the precomputation cost in Section~\ref{subsec:precompute}. However this cost is
unimportant because we only need to do it once for the chosen $p$ for
a given kernel. In the evaluation stage, step a) has the usual
$\Ocal(N)$ complexity of KIFMM, and takes slightly more time compared
to free-space KIFMM, because the near-field evaluation requires S2T
direct summation of interactions across the periodic boundaries. This
extra cost scales as $\Ocal(sN_B)$, where $s$ is the maximum number of
points allowed in one leaf box in the adaptive octree, and $N_B$ is
the number of points located in the leaf boxes adjacent to the
periodic boundary in the octree. The value of $N_B$ depends on the distribution of
source and target points. Step b) takes negligible time because
$\phi^{B_0,u}$ is simply an array of $\Ocal(p^2)$ elements. Step c) is
a simple matrix-vector multiplication with a $\Ocal(p^4)$ cost, which
is fast as usually $p<20$. Step d) has a
$\Ocal(T p^2 )$ cost where $T$ is the number of target points in
$B_0$, but we found this cost still much smaller than the cost of step
a). Timing results are reported in
Section~\ref{subsec:stokestiming}.

Algorithm~\ref{alg2} is straightforward to implement with
multi-threading utilities like \verb|OpenMP|, because the pattern of
step c) and step d) are `embarrassingly parallel', requiring no
communication. In the case of a distributed memory cluster with
\verb|MPI|, no modification to the algorithm is necessary. In fact, all
inter-node communications are handled by the underlying KIFMM package
in step a). Step c) has little cost and is easy to calculate on every
node, and in step d) every node only needs to calculate a subset of
$q_{\Fcal}^t$ located on its own memory.

The error of the periodizing method observes the same theoretical
error estimates found for the original KIFMM
work \cite{ying_kernel-independent_2004} because we designed the
algorithm to reuse the data calculated by KIFMM, and the calculation
of $\MLOP$ uses the same KIFMM formulation. Also, the $K^{P,F}$ is
straightforward to evaluate to machine precision by the Ewald method
with its well-known error analysis \cite{lindbo_spectrally_2010}, or by
direct summation, depending on the geometry. There is no need to
repeat the error analysis in this work, and we show in
Section~\ref{sec:results} that we achieved the same error as the
underlying PVFMM package.

In the following we show accuracy and timing results for the Laplace
kernel and the Stokes velocity kernel\footnote{It is sometimes
  referred to as the Stokes single layer kernel in the boundary
  integral method.}:
\begin{align}\label{eq:kexamples}
  K&=\frac{1}{\ABS{\bx_{t}-\by_{s}}},
  \\ \bK&=\frac{1}{8\pi}\left(\frac{\bI}{\ABS{\bx_{t}-\by_{s}}}
  + \frac{(\bx_{t}-\by_{s})(\bx_{t}-\by_{s})^T}
  {\ABS{\bx_{t}-\by_{s}}^3}\right).
\end{align}

\section{Results}
\label{sec:results}
\subsection{The backward stable solver for $\MLOP$.}
\label{subsec:SVD}

The matrix $\Acal$ in Eq.~(\ref{eq:MLOPAxb}) is well-known to be
numerically rank-deficient with increasing $p$ in the method of
fundamental solutions (MFS) \cite{barnett_stability_2008}. Table~\ref{tab:cond} shows the condition number $\kappa$ and numerical rank of $\Acal$ for the Laplace and
Stokes velocity kernels, where $r$ denotes the numerical rank of $\Acal$, and
$\dim\Acal\propto6(p-1)^2+2$ represents the dimension of the square
matrix $\Acal$. The check and equivalent points are located on a cube
with edge length $1.05$ and $2.95$, respectively.
\begin{table}[ht]
	\centering
	\caption{Conditioning  of the matrix $\Acal$. The condition number grows approximately exponentially with increasing $p$, as is well known in the MFS (Method-of-Fundamental-Solutions) community \cite{barnett_stability_2008}.   We identify singular values less than $\epsilon s_{max} \dim\Acal$ as indicating rank deficiency, where $s_{max}$ is the largest singular value and $\dim\Acal$ is the matrix dimension. $\epsilon=2^{-52}$ is the floating point relative accuracy, which is defined as the distance from $1.0$ to the next larger double precision number. This threshold is the default setting of both \texttt{numpy} and \texttt{MATLAB}, and also appears in \cite{press_numerical_2007}. } 
	\begin{tabular}{cc|cccccc}
		\hline 
	&	$p$ & 6 & 8 & 10 & 12 & 14 & 16 \\ \hline
Stokes&	$ \kappa$ & $2.34E10$ & $1.26E13$ & $4.35E16$ & $3.57E19$ & $3.60E20$ & $1.39E21$ \\
	&	$r/\dim\Acal$ & $456/456$ & $881/888$ & $1320/1464$ & $1609/2184$ & $1822/3048$ &  $1977/4056$ \\
Laplace& $ \kappa$ & $1.13E9$ & $3.61E12$ & $1.19E16$ & $4.20E19$ & $3.64E19$ & $1.75E20$                             \\
	&	$r/\dim \Acal$ & $152/152$ & $296/296$ & $444/488$ & $550/728$ & $628/1016$ & $688/ 1352$									\\
\hline
	\end{tabular}
        \label{tab:cond}
\end{table}
The rank-deficiency problem was initially ameliorated with
Tikhonov regularization \cite{ying_kernel-independent_2004}. In
PVFMM \cite{malhotra_pvfmm_2015}, a backward stable solver is used to
improve the accuracy from $\sim 10^{-9}$ to $\sim 10^{-14}$ without
regularization.

The backward stable solver is well-known to the numerical analysis
community, but perhaps less so for more general readers.  We find that
backward stability is crucial, otherwise the accuracy may stagnate at
$10^{-7}$ (or worse due to the large condition number of
$\Acal$. To assist in implementing the method proposed in this paper,
we include here a brief description of the backward stable solver, and
in \ref{sec:svdcompare} a brief comparison of open-sourced
implementations.  The solver works as follows for Eq.~(\ref{eq:MLOPAxb}):
\begin{enumerate}[1.]
	\item Construct $\Acal$, and calculate its SVD: $\Acal=\bU\bS\bV^T$.
	\item Calculate the approximate pseudoinverse $\bS_\epsilon^+$
          of $\bS$, by inverting all non-zero singular values in $\bS$
          whose absolute value is greater than $\epsilon_{SVD}$, and setting
          all others to zero.
	\item Calculate $\bphi_{j}^{B_0,d}$ for each $j$ in the order
          implied by the parenthetical nesting: $\bphi_{j}^{B_0,d} =
          \left( \bV\bS_\epsilon^+ \right)\left(\bU^T \bq_j^{B_0,d}
          \right) $, or $\bphi_{j}^{B_0,d} = \bV\left(
          \left(\bS_\epsilon^+\bU^T \right) \bq_j^{B_0,d} \right) $.
\end{enumerate}
The main point is to avoid the explicit construction of the
approximate pseudoinverse matrix
$\Acal_{\epsilon}^+=\bV\bS_\epsilon^+\bU^T$ which yields a significant loss of accuracy.

There is a simple explanation as to why high accuracy can be achieved despite the numerical rank deficiency of matrix $\Acal$. Equation~(\ref{eq:MLOPAxb}) means the density $\phibb_{j}^{B_0,d}$ should be determined to match the condition on the check surface, and an almost singular $\Acal$ means the solution is (numerically) not unique. However, to generate the field we need, there is no need to distinguish between possible solutions. As long as Eq.~(\ref{eq:MLOPAxb}) is accurately satisfied, the $\phibb_{j}^{B_0,d}$ does equally well for the ensuing potential calculations. The backward stability of the solver guarantees that any
solution given by the solver satisfies Eq.~(\ref{eq:MLOPAxb}) to machine precision, and therefore any solution will do. In this sense, the operator $\MLOP$ is not, nor need be, unique. As a consequence, there is no need to use the same solver in the precomputing of $\MLOP$ and
the actual KIFMM calculation, as long as they are both backward-stable. For example, one can solve for $\MLOP$ conveniently in MATLAB with the backslash operator and safely use the results with PVFMM.
 
\subsection{Laplace Kernel: the Madelung constant for a 1D crystal}
\label{subsec:laplace1D}

Here we report accuracy test results for computing the 1D Madelung constant. This is a simple example but we include more details on our implementation of Algorithm~\ref{alg2}.

First, the compatibility condition should be considered. For a general singly periodic Laplace problem, the compatibility condition is simply the neutrality condition in Eq.~(\ref{eq:compatcondition}). The second step is to solve Eq.~(\ref{eq:M2Lsolve}). Usually this step is done by first calculating the periodic Ewald kernel $K^{E}$ and subtracting the $\Ncal^{B_0}$ contributions. However, $K^E$ for singly periodic Laplace geometry involves an incomplete modified Bessel function of the second kind \cite{tornberg_ewald_2015} $L(u,v)=\int_{1}^{\infty}\exp{(-ut-v/t)}/tdt$, which is very hard to evaluate uniformly to machine precision for any $u,v$, though there are attempts in the literature \cite{harris_incomplete_2008,alford_calculation_2005}.

Instead of invoking Ewald summation, we note that due to the neutrality condition the naive direct sum of $n$ periodic images in singly periodic geometry converges as $1/n^2$. By direct summation we can reach an error $\sim10^{-10}$ with $n=10^5$, which poses no problem on modern hardware.\footnote{Stokes velocity kernel observes the same convergence rate. Therefore for singly periodic problems we do not need to deal with the special function $L(u,v)$ if precomputation time is not a concern.} The periodic kernel $K^{P,F}$ is then converted to a direct summation:
\begin{align}\label{eq:KPFdirect}
K^{P,F} = \sum_{\ell< p_z < n} \frac{1}{\ABS{\bx_{t}-\by_{s} - \bp}} + \frac{1}{\ABS{\bx_{t}-\by_{s} + \bp}} 
\end{align} 
where we assume the periodic condition is imposed in the $z$ direction, $\bp=(0,0,p_z)$, and the two terms with $\pm\bp$ means the image at positive and negative z directions. It is worth noting that for a given target/source pair ($\bx_{t},\by_{s}$), $K^{P,F}$ scales as the harmonic series $\sum_{j=l}^n \frac{1}{j}$, and is slowly divergent as $\log n$.  
When we attach this far field to the near field, we effectively approximate the true periodic system with a finite chain of $n$ images on both sides of $B_0$. Since a real system satisfying charge neutrality condition converges as $1/n^2$, we can use $n=10^6$ to get convergent results and the divergent $K^{P,F}$ does not matter. In general when $K^{P,F}\gg 1$, truncation error may occur due to finite floating point precision. However, this is not a problem in this case because when $n=10^6$, $K^{P,F}$ is only on the order of $10$. $\MLOP$ is found with the backward stable solver for different
$p\in[6,16]$.  In this example, we demonstrate the accuracy and convergence of Algorithm~\ref{alg2} by calculating the Madelung constant $M_{1D}$ for a 1D crystal in 3D space. It is defined as the potential generated by a periodic chain of equispaced alternating
$\pm1$ unit charges on one charge in the chain. Mathematically it is given as a series summation: 
\begin{align}
M_{1D} = \sum_{i\in Z/0} \frac{(-1)^i}{i} = -2\ln2 = - 1.386 294 361 119 890 \cdots 
\end{align}
We place 4 charges in the unit cubic cell and assume the unit cubic cell is periodic in $z$ direction: $\phi^s=\{1,-1,1,-1\}$, $\by_{s}=\{(0.5,0.5,0.125),(0.5,0.5,0.375),(0.5,0.5,0.625),
(0.5,0.5,0.875) \}$, and we evaluate the potential at the two positive charges $q_1$ and $q_3$: $\bx_t=\{(0.5,0.5,0.125), (0.5,0.5,0.625) \}$. In this case, the analytical solution at $\bx_t$ is $-8\ln2=4M_{1D}$. 

\begin{figure}[h]
\centering
\includegraphics[scale=1.0]{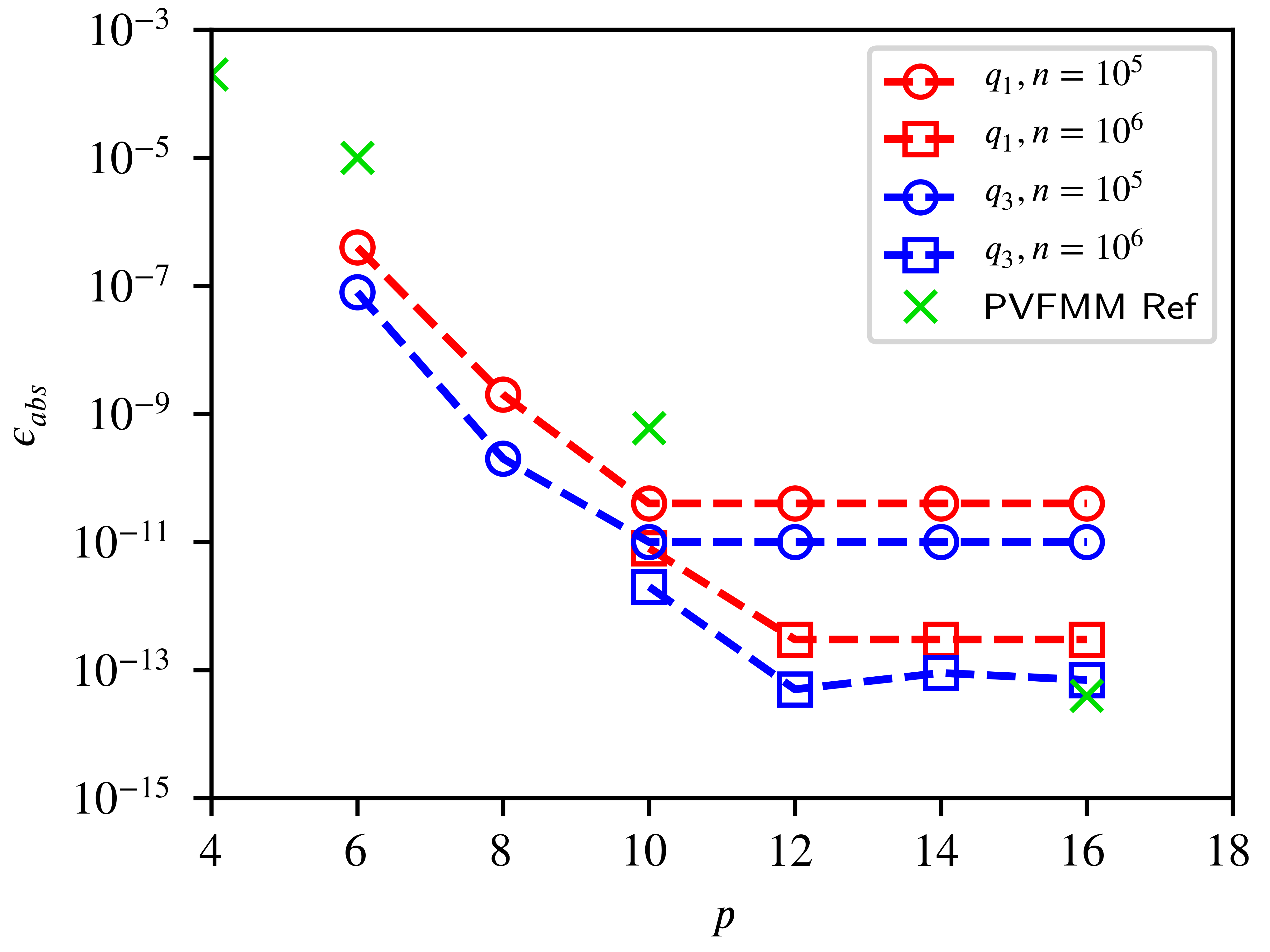}
\caption{Absolute error for the potentials $q_1$ and $q_3$ evaluated for the two positive charges, compared to the analytical solution $M_{1D}$. \textsf{PVFMM Ref} refers to the error achieved by PVFMM package at the same $p$  \cite{malhotra_pvfmm_2015}, although PVFMM is not used in this calculation. For $n=10^5$,    $\epsilon_{abs}$ stagnates at $\sim10^{-11}$, while for $n=10^6$,   close to machine precision is achieved. }
\label{fig:Laplace1D}
\end{figure}

Since it is a tiny problem of only 4 points, there is no need to invoke FMM. $q_{\Ncal}^t$ is evaluated as a direct sum and $\phi^{B_0,u}$ is evaluated separately. The accuracy results are
reported as the absolute error $\epsilon_{abs} = \ABS{q^t+8\ln2}$ in Figure~\ref{fig:Laplace1D}, where $n$ is the number of periodic images directly summed in $K^{P,F}$. In all tests the splitting between near and far-field is chosen at $\ell=2$ layers of image boxes. \textsf{PVFMM Ref} refers to the error achieved by PVFMM package at the same $p$  \cite{malhotra_pvfmm_2015}, although PVFMM is not used in this calculation.
For $n=10^6$, close to machine precision is achieved. There is a systematic error for different $n$, where $q_3-q_1$ stagnates beyond some $p$. This is because we constructed the $\MLOP$ with a finite $n$, which approximate the true periodic system with a finite chain. Due to the asymmetric charge distribution in $B_0$, which is related to the net dipole moment of $B_0$, $q_1$ and $q_3$ do not feel exactly the same environment in this finite chain of $2n+1$ periodic boxes. We can calculate that $q_3-q_1 \approx  - \frac{1}{2n^2}$, neglecting the $4\pi$ factor for Laplace kernel. However, since fixing this error does not improve the order of magnitude of accuracy, we neglect it.

All calculations for $\MLOP$ with $n=10^5$ complete within $\sim 15$
minutes on an Intel Xeon E5-2697 v3 with 14 cores at 2.60GHz. For
single precision accuracy, $p=6$ is enough and in this case the total
number of equivalent points is only $6(p-1)^2+2=152$. Accuracy close
to the double precision limit can be achieved with large $p$ and
$n=10^6$, and the calculation takes less than 2 hours.
$\MLOP$ takes this long is a special case due to the inapplicability of 1D Ewald formula. In general the computation for $\MLOP$ completes within $10\sim 100$ seconds. We discuss the precomputation cost in more detail in Section~\ref{subsec:precompute}.

\subsection{Stokes Kernel: doubly periodic geometry in 3D space}
\label{subsec:stokes2P}
In this section, we present results for 3D Stokes flows, using the
velocity kernel in Eq.~(\ref{eq:kexamples}), in $B_0$ with periodicity in
the $x-y$ plane. Here, the compatibility condition takes the
neutrality form if we require the fluid velocity $\bu$ to {approach}
finite values as $z\to\pm\infty$ \cite{lindbo_fast_2011}. More
precisely, the net force in $B_0$ must be zero. In this case the
simple direct sum in Eq.~(\ref{eq:KPFdirect}) is not tenable as the planar
periodic summation of force-neutral cells converges very slowly,
sometimes as slow as $1/\sqrt{n}$ \cite{borwein_convergence_1985},
where $n$ is the number of cells summed. In this case, a doubly periodic
Ewald sum must be invoked to calculate $\MLOP$. We skip the details
about implementing a pairwise Ewald sum, which has been thoroughly
discussed in the literature \cite{lindbo_fast_2011}.

With the pairwise Ewald formulation, we have:
\begin{align}\label{eq:KPFEwald}
\bK^{P,F} = \bK^{Ewald} - \sum_{-\ell\leq p_z \leq \ell } \bK(\bx_{t},\by_{s}+\bp),
\end{align} 
where the splitting between $\Ncal^{B_0}$ and $\Fcal^{B_0}$ is still
chosen as $\ell=2$. $\bK^{Ewald}$ is calculated by splitting the sum
into real and wave space by the splitting parameter $\xi$, and then
the real and wave space sums are both truncated at the desired
accuracy.  In our calculation of $\MLOP$ the splitting parameter $\xi$
and truncation parameters are taken such that the accuracy is close to
machine precision, independent of the Ewald parameters. Several
methods are available to accelerate the Ewald evaluation using the FFT
 \cite{lindbo_fast_2012}, but we follow a simple pairwise scheme
because we need only evaluate Ewald sums from the equivalent surfaces
to check surfaces in Eq.~(\ref{eq:M2Lsolve}), and only in the precomputing
stage to solve $\MLOP$. $\MLOP$ is determined for $p\in[6,16]$.

A physical consequence of the compatibility condition is that the net
flux of flow across the periodic $x-y$ plane must be zero:
\begin{align}\label{eq:fluxz}
\int_{B_0} u_z(\bx) d^3\bx = 0.
\end{align}
Therefore, this can be used as a measurement of numerical error. However, a direct numerical evaluation of $u_z^{net}$ is non-trivial because $\bu(\bx)$ is not smooth in $B_0$. If a point force $i$ is placed at $\by_i$, $\bu(\by_i)$ is singular. To circumvent this singularity problem, we calculate the following surface integral as a physically equivalent measurement of error:
\begin{align}\label{eq:zflux}
			u_z^{net}&=\int_{S_z} u_z(\bx) dS_z = 0, \quad S_z=\{(x,y,z) | x\in[0,1],y\in[0,1],z=0\}
\end{align}
The surface $S$ is chosen to be at $z=0$, away from all point forces inside $B_0$, so that $u_z(\bx)$ is smooth on $S$. $u_z^{net}$ therefore can be accurately integrated with a simple quadrature rule, and can be used as an error measurement.

\begin{figure}[h]
\centering
\includegraphics[width=0.45\linewidth]{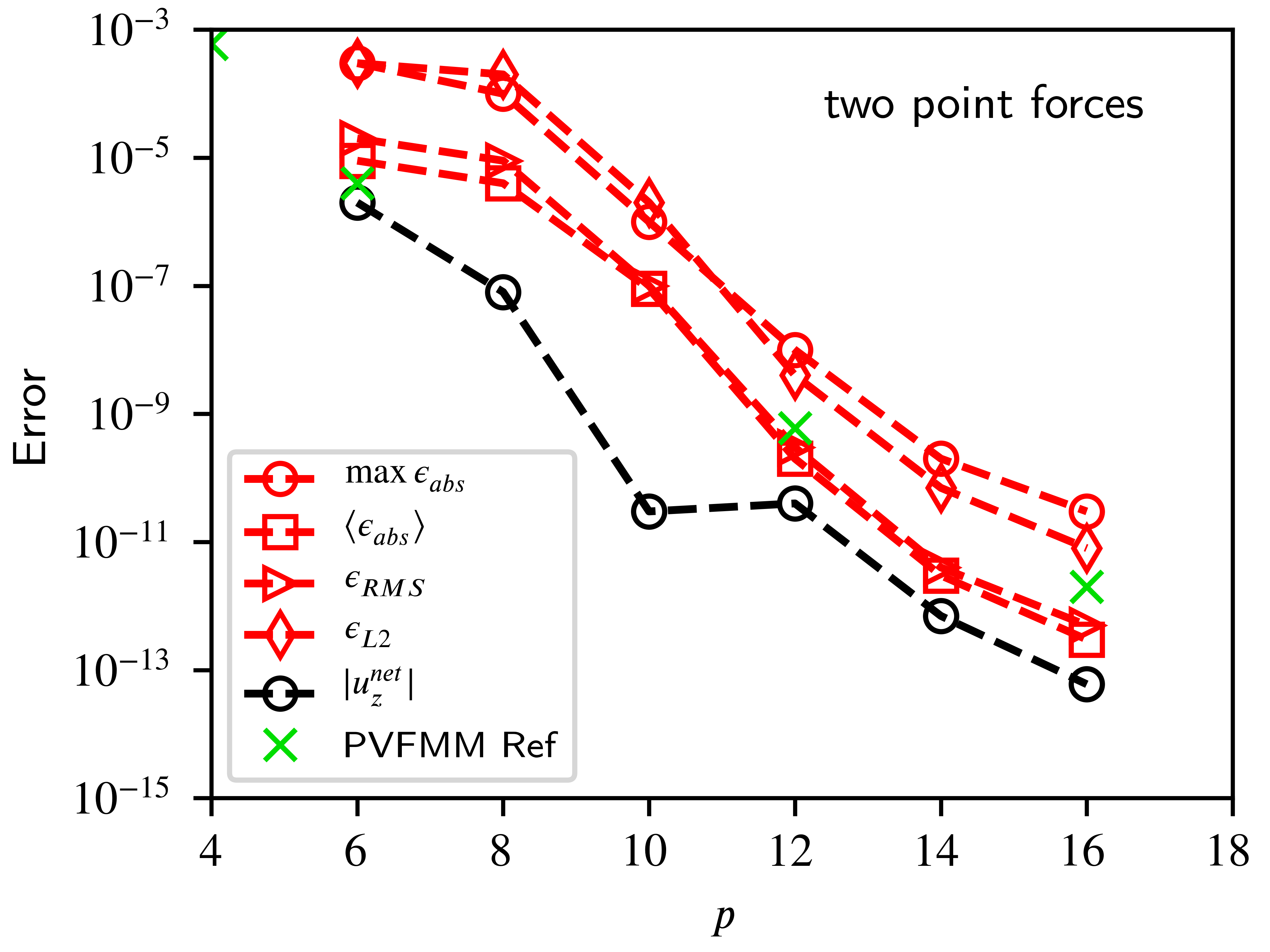}
\caption{Accuracy test results for two point forces placed
  in $B_0$, for induced flows
  periodic in the $x-y$ plane. The red symbols are accuracy results
  with different error measurements, the green symbols are the reference
  error achieved by PVFMM \cite{malhotra_pvfmm_2015}. The black symbols
  are for the net flow in the $z$ direction.
}
\label{fig:Dipole2P}
\end{figure}

To assess accuracy we place two unit point forces $\bff_1,\bff_2$
in $B_0$, with $\bff_1+\bff_2=0$
  imposed because the convergence in doubly periodic geometry requires
  a zero net forcing in $B_0$. We fix $\bx_1=(0.7,0.6,0.5)$, $\bff_1=(1,2,3)/\sqrt{14}$, $\bx_2=(0.2,0.8,0.7)$, $\bff_2=-\bff_1$ so that all the following accuracy tests are generated with the same forcing. This particular choice of point force has no effect on the accuracy test results. $\bff_1$ and $\bff_2$ generate the
  velocity field of a force dipole in the far field of $B_0$. The
  induced flow velocity $\bu$ is evaluated at a set of target points
  $\bx_{ijk}=(c_i,c_j,c_k)\in B_0$. The source points are the
  Clenshaw-Curtis (Chebyshev) quadrature points mapped to the range
  $[0,1]$. As illustrated in Figure~\ref{fig:points}, $97$ points are
  placed in each dimension with the total number of target points is
  $97^3=912673$. The surface integration of net flow is approximated
  as a discretized tensorial sum, where $w_i,w_j$ are the numerical
  integration weights:
\begin{align}\label{eq:uint}
u_z^{net} = \int_{S_z} u_z(\bx) dS_z \approx \sum_{i,j} \bu(\bx_{ij,k=0}) w_iw_j
\end{align}

Errors are calculated by comparing $\bu(\bx_{ijk})$ with the
`accurate' result calculated by the direct pairwise Ewald sum. To ease the
comparison with results in the literature, we report errors based on
four different definitions:
\begin{align}
	&\text{Maximum absolute error: } &\max\epsilon_{abs} =& \max_{i\in \dim \bU} \ABS{U_i-U_i^{Ewald}}, \\
	&\text{Average absolute error: } &\AVE{\epsilon_{abs}} =& \frac{1}{\dim \bU} \sum_i\ABS{U_i-U_i^{Ewald}}, \\
	&\text{Root mean square error: } &\epsilon_{RMS} =& \sqrt{\frac{1}{{\dim \bU}}\sum_i(U_i-U_i^{Ewald})^2}, \\
	&\text{Relative L2 error: } &\epsilon_{L2} =& \lVert{\bU-\bU^{Ewald}}\rVert_2 /\lVert{\bU^{Ewald}}\rVert_2,
\end{align}
where $\bU$ is the large vector consists of all $\bu_{ijk}$, i.e.,
$\bU=(\bu_{000},\bu_{001},\cdots)$. With this definition, the above
errors are actually defined in the sense of \emph{per velocity
  component}. The results are reported in Figure~\ref{fig:Dipole2P}.
These results show the same level of accuracy as the PVFMM package
used to evaluate $q_{\Ncal}^t$. In comparison, the Spectral Ewald
method \cite{lindbo_spectrally_2010} can achieve $\epsilon_{RMS} \to
10^{-14}$, and our method can be very close in accuracy without the large
cost of the 3D FFT mesh. The net flow condition $u_z^{net}=0$ is also
properly imposed.

\subsection{The Stokes kernel: unit point force in a triply periodic geometry}
\label{subsec:stokes3P}

As in the last section, we start from the compatibility condition. For
triply periodic Stokes flow the sum of all point forces need not be
zero in the periodic box. The net forcing is balanced by a pressure
gradient governed by the Stokes equation, which is usually
\emph{chosen} such that the net flow of material in space is zero. The
physical explanation behind this zero net flow condition is that, the
numerical solution is describing an experiment in an infinitely large
fixed container, where no material flows in or out of it. A more
detailed discussion about convergence can be found in \cite{obrien_method_1979}. The zero net flow condition
is imposed by setting the constant $C=0$ in Eq.~(\ref{eq:compatcondition}):
\begin{align}\label{eq:flowvolume}
\int_{B_0} \bu(\bx) d^3\bx = 0.
\end{align}
To guarantee this condition, especially for the case of non-zero net forcing, an Ewald sum must be used where the above zero-flow condition is guaranteed by removing the $\bk=0$ term in the wave space sum. In this case, $\bK^{P,F}$ is still calculated as Eq.~(\ref{eq:KPFEwald}) where $\ell=2$ is fixed. We follow the well-known Hasimoto convention of $\bK^{Ewald}$ \cite{hasimoto_periodic_1959,lindbo_spectrally_2010}.

To assess the accuracy, we avoid the numerical integration of Eq.~(\ref{eq:flowvolume}) due to the singularities of $\bu(\bx)$ as discussed in the previous subsection. Instead, we numerically integrate the following three fluxes in the $x,y,z$ directions on the $x=0,y=0,z=0$ surfaces of $B_0$, respectively:
\begin{align}
	u_x^{net}&=\int_{S_x} u_x(\bx) dS_x = 0, \quad S_x=\{(x,y,z) | x=0,y\in[0,1],z\in[0,1]\}, \\
	u_y^{net}&=\int_{S_y} u_y(\bx) dS_y = 0, \quad S_y=\{(x,y,z) | x\in[0,1],y=0,z\in[0,1]\}, \\
		u_z^{net}&=\int_{S_z} u_z(\bx) dS_z = 0, \quad S_z=\{(x,y,z) | x\in[0,1],y\in[0,1],z=0\}.
\end{align}
The same Chebyshev quadrature as in Eq.~(\ref{eq:uint}) is used to calculate $u_x^{net},u_y^{net}$ and $u_z^{net}$ numerically. 

\begin{figure}
\centering
\includegraphics[width=0.45\linewidth]{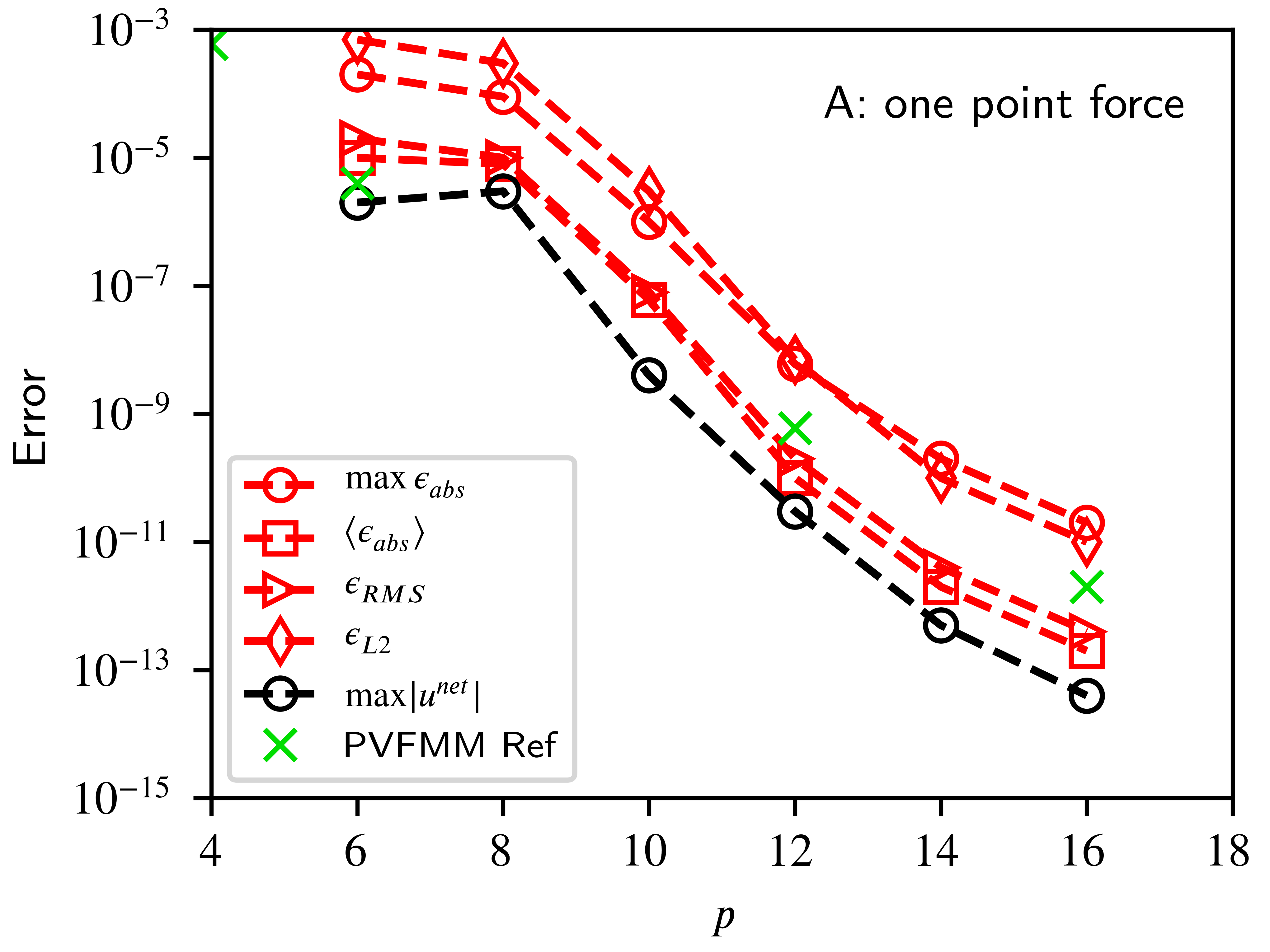}
\includegraphics[width=0.45\linewidth]{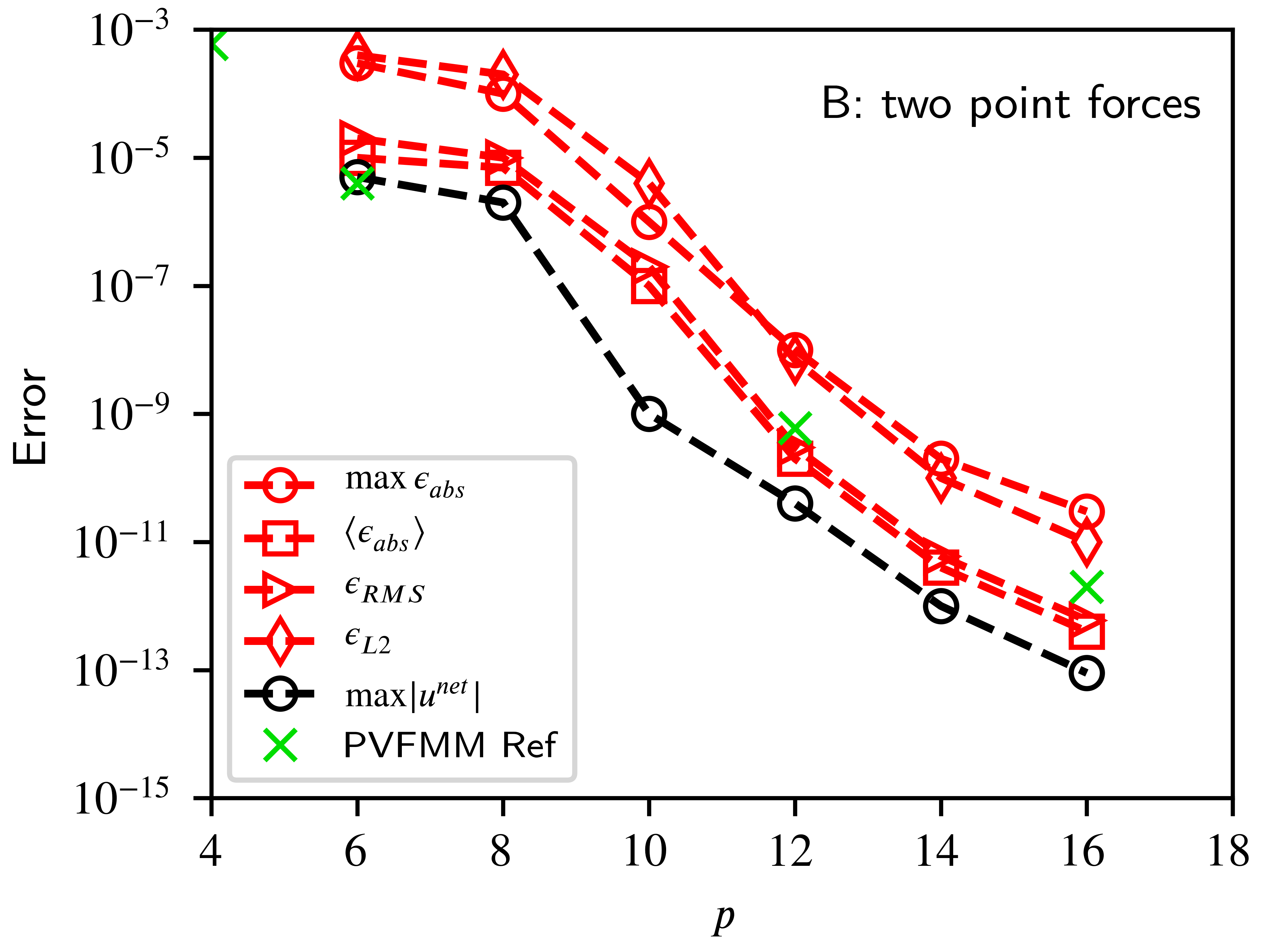}
\caption{Accuracy test results for the triply periodic boundary condition test case. The red symbols are
  accuracy results using different error measurements, the green symbols
  are the reference error achieved by
  PVFMM \cite{malhotra_pvfmm_2015}. $\max\ABS{u^{net}}$ is the maximal
  component of $\bu^{net}$.}
\label{fig:force3P}
\end{figure}

We present accuracy test results for two cases, where the target points are the same set of points as in the last subsection. The triply periodic boundary condition allows a net force in $B_0$. Therefore, in addition to the test of two point forces shown before, we added a test with $\bx=(0.7,0.6,0.4)$, $\bff=(1,2,3)/\sqrt{14}$. Again, this particular choice of point force is chosen only to ensure the same input of all test cases to help the comparison later with different $\ell$. This choice has no impact on the accuracy test results. Figure~\ref{fig:force3P} A shows the accuracy test results for the point force test. Figure~\ref{fig:force3P} B shows the accuracy test results for the test of two point forces, as specified in the last subsection. The results show the same level of accuracy as in the last subsection.

\subsection{The Stokes kernel: Timings for large systems }
\label{subsec:stokestiming}

\begin{figure}
	\centering
	\includegraphics[scale=1.0]{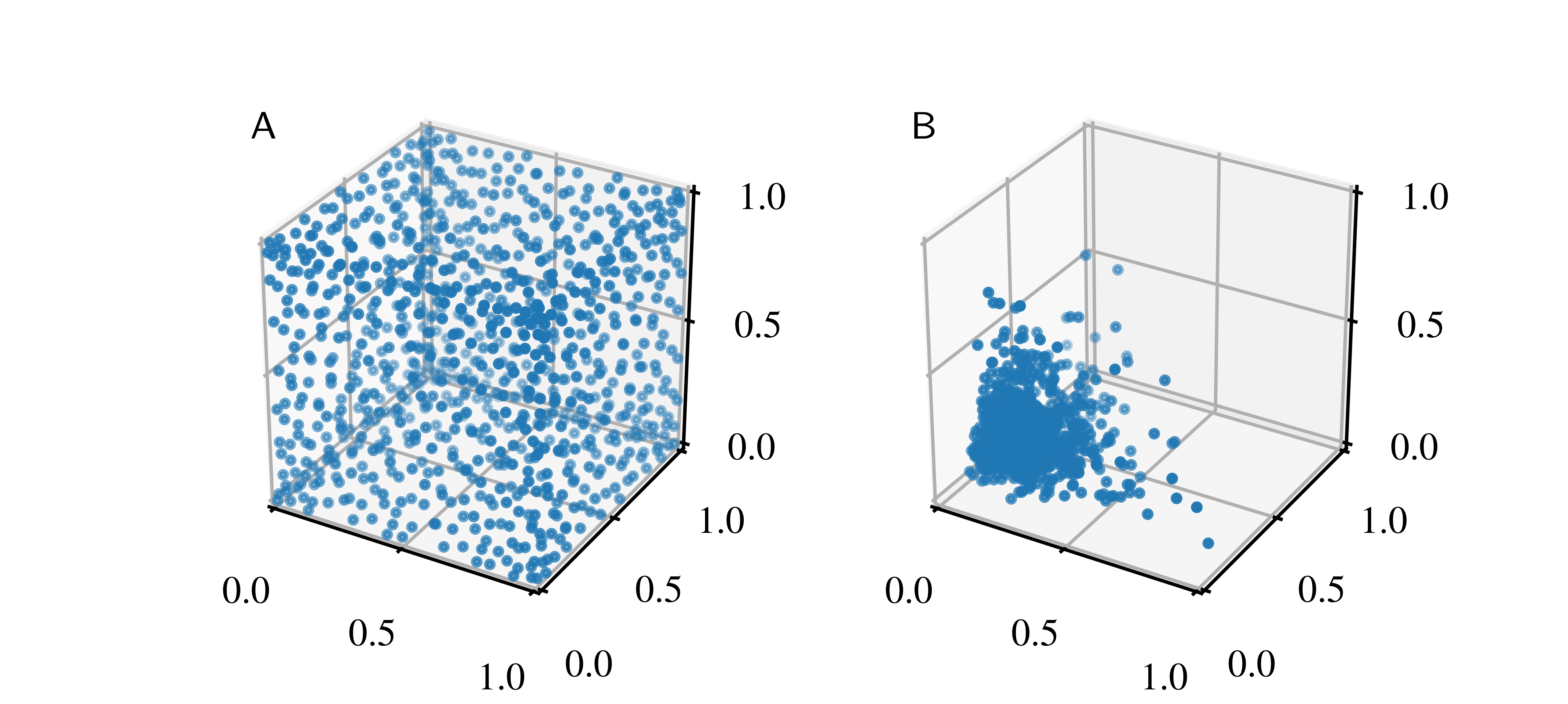}
	\caption{Different distributions of source and target points.  A:
          Chebyshev points $\bx_{ijk}=(c_i,c_j,c_k)$. B: Random points
          $\bx_{ijk}$, where $(x_i,x_j,x_k)$ are generated from a
          lognormal random distribution with standard
          parameters $(0.2,0.5)$. Only $0.1\%$ of the points in the
          computations are shown in both Figure A and B.}
	\label{fig:points}
\end{figure}

We now provide timing results for steps a), c), and d) in
Algorithm~\ref{alg2} for uniform and non-uniform distributions of
source and target points. In this part, the same
set of points are used for both the source and target points. The
uniform points are chosen to be the same set of $97^3$ points
$\bx_{ijk}$ described previously, as shown in Figure~\ref{fig:points}
A. The non-uniform point distributions are a set of $97^3$ points with each
component of $\bx_{ijk}$ generated from a lognormal distribution with
standard parameters $(0.2,0.5)$, as shown in
Figure~\ref{fig:points} B.  A random point force $\bff$ is placed on
each source point $\bx_{ijk}$, with each component
generated from an uniform distribution on $[-0.5,0.5)$. For tests with
  singly and doubly periodic boundary conditions, a uniform force is
  added to each point force to guarantee the force neutrality
  condition.

\begin{table}[h]

	\centering
	\caption{Timing (in seconds) results for $\sim10^6$ Chebyshev
          points in Figure~\ref{fig:points} A. SP, DP, and TP refer to
          singly, doubly and triply periodic boundary
          conditions. $\%_{cost}$ refers to the cost of the far-field
          calculation relative to the near-field calculations, and is
          calculated as: $\%_{cost} =
          \tau_{\Fcal}/\left(\tau_{tree}+\tau_{\Ncal}\right)$. $\tau_{tree}$
          is independent of boundary condition. In reality, $\tau_{M2L}\sim 10^{-4}$s for $p=6$ and $\sim 10^{-2}$s for $p=16$. It is much less than other parts, and so is not included in this table and is ignored in $\%_{cost}$ calculations.}
	\begin{tabular}{c|cc|cccc|cccc|cccc}
		\hline
		    &\multicolumn{2}{c|}{ {Free-Space}} & \multicolumn{4}{c|}{ {SP}} & \multicolumn{4}{c|}{ {DP}} & \multicolumn{4}{c}{ {TP}} \\
		$p$ & $\tau_{tree}$ & $\tau_{\Ncal}$ & $\tau_{tree}$  & $\tau_{\Ncal}$ & $\tau_{\Fcal}$ & $\%_{cost}$& $\tau_{tree}$ &   $\tau_{\Ncal}$ & $\tau_{\Fcal}$ & $\%_{cost}$ &  $\tau_{tree}$&  $\tau_{\Ncal}$ & $\tau_{\Fcal}$ & $\%_{cost}$\\
		\hline
		$6$ & 0.44 & 2.22 & 0.44 & 2.29 &0.18 & 7\% & 0.42 &2.49 &0.22 & 7\%& 0.44& 2.92 & 0.14 & 4\% \\
		$8$ & 0.45 & 3.33  & 0.48 &3.48 &0.26 & 7\% & 0.45 &3.80 &0.26 & 6\% & 0.45&  4.32 & 0.26 & 5\% \\
		$10$ & 0.55 & 4.31 & 0.52&4.83 &0.44 & 8\% & 0.55 & 5.23 &0.59 & 10\% &0.55 &  5.67 & 0.42 & 7\% \\
		$12$ & 0.66 & 7.89 & 0.62 &8.47 &0.70 & 8\% & 0.62 & 8.71  &0.68 &7\% &0.62 & 9.76 & 0.70 & 7\% \\
		$14$ & 0.83 & 13.31 & 0.86 &14.84 &0.95 & 6\%& 0.83  &16.01  &0.97 &6\%  &0.85 & 16.67 & 0.96 & 5\%\\
		$16$ & 0.98 & 26.32& 1.14 &28.91 &1.21 & 4\%& 0.97 &31.17  &1.23 &4\% & 0.98& 32.60 & 1.24 & 4\%\\
		\hline
	\end{tabular}
\label{tab:StokesTiming}
\end{table}

\begin{table}[h]

	\centering
	\caption{Timing (in seconds) results for $\sim10^6$ non-uniform random points in Figure~\ref{fig:points} B. SP, DP, TP and $\%_{cost}$ are the same as Table~\ref{tab:StokesTiming}. In reality, $\tau_{M2L}\sim 10^{-4}$s for $p=6$ and $\sim 10^{-2}$s for $p=16$. It is much less than other parts, and so is not included in this table and is ignored in $\%_{cost}$ calculations. Also in this case the near field of TP takes less time to compute than DP and SP is a special case depending on the point configuration. In general we should expect TP to be slower than DP and SP. } 
	\begin{tabular}{c|cc|cccc|cccc|cccc}
	\hline
	&\multicolumn{2}{c|}{ {Free-Space}} & \multicolumn{4}{c|}{ {SP}} & \multicolumn{4}{c|}{ {DP}} & \multicolumn{4}{c}{ {TP}} \\
	$p$ & $\tau_{tree}$ & $\tau_{\Ncal}$ & $\tau_{tree}$  & $\tau_{\Ncal}$ & $\tau_{\Fcal}$ & $\%_{cost}$& $\tau_{tree}$ &   $\tau_{\Ncal}$ & $\tau_{\Fcal}$ & $\%_{cost}$ &  $\tau_{tree}$&  $\tau_{\Ncal}$ & $\tau_{\Fcal}$ & $\%_{cost}$\\
	\hline
	$6$ & 0.47 & 2.74 & 0.47  & 2.67 & 0.14 & 4\% & 0.47 & 2.73 & 0.14& 4\%&0.45 & 2.68 &0.17 &5\% \\
	$8$  & 0.51 & 4.25 & 0.49  & 4.23& 0.32& 7\% & 0.50 & 4.53 & 0.28& 6\%& 0.53& 4.53&0.27 & 5\% \\
	$10$  & 0.68 & 5.55 &  0.61 &5.8 & 0.49& 8\% & 0.63 & 6.00 & 0.49& 7\%& 0.68& 5.74 &0.48 & 8\% \\
	$12$  & 0.73 & 11.75 & 0.68  &10.37 &0.67 & 6\% & 0.68 & 10.40 &0.80 & 7\%&0.71 &10.15  &0.68 & 6\% \\
	$14$  & 0.90 & 17.86 & 0.87  &18.27 &0.96 & 5\% & 0.83 & 17.65 &0.95 & 5\%&0.87 &17.51  &1.02 & 6\% \\
	$16$  & 1.03 & 33.09 & 1.05  &35.49 &1.24 & 3\% & 1.00 & 36.13 &1.23 & 3\%&1.07 &34.67  &1.23 & 3\% \\
	\hline
\end{tabular}
	\label{tab:StokesTimingRandom}
\end{table}

The results reported in Table~\ref{tab:StokesTiming} and
Table~\ref{tab:StokesTimingRandom} are for uniform and non-uniform cases,
respectively. Timings are reported as seconds of wall clock time,
measured on a machine with two sockets of Intel(R) Xeon(R) CPU E5-2643
v3 @ 6 core 3.40GHz and 128GB memory. \verb|OpenMP| is enabled but
\verb|Hyper-Threading| is disabled to give more consistent timing
measurements. The maximum number of points in a leaf octree box $s$ is
set to $s=1000$ in the PVFMM routine, and the splitting between
$\Ncal^{B_0}$ and $\Fcal^{B_0}$ is fixed at $\ell=2$. $\tau_{tree}$
and $\tau_{\Ncal}$ are the FMM tree construction time and the FMM
evaluation time for $\Ncal^{B_0}$ in step a), respectively. $\tau_{M2L}$ is the time
to calculate $\phi^{B_0,d}(\by) = \MLOP \phi^{B_0,u}$ of step
c). $\tau_{\Fcal}$ is the time to calculate $q_{\Fcal}^t$ and add it
to $q_{\Ncal^t}$ in step d). $\%_{cost}$ is defined as the relative
cost of the far-field evaluation: $\%_{cost} =
\left(\tau_{M2L}+\tau_{\Fcal}\right)/\left(\tau_{tree}+\tau_{\Ncal}\right)$. In
reality, $\tau_{M2L}\sim 10^{-4}$ seconds and is much less than other
parts, and so is not included in Table~\ref{tab:StokesTiming}
and is ignored in $\%_{cost}$ calculations.


The time $\tau_{tree}$ is independent of boundary condition because the octree
is `copy-and-pasted' for all boxes in $\Ncal^{B_0}$. Such
`copy-and-paste' also saves a significant amount of time in near-field
evaluations because only the S2T interactions (direct interactions)
from the leaf boxes across the periodic boundaries require additional
cost compared to the free-space boundary conditions. Other non-leaf
octree boxes can be evaluated at the precomputing stage of PVFMM for
the octree data\footnote{This is different from the precomputing stage  of $\MLOP$ in Algorithm~\ref{alg2}, but is also independent of the   source and target locations in $B_0$. Such precomputing techniques are   also implemented in other KIFMM packages such as   ScalFMM \cite{bramas_optimization_2016,agullo_task-based_2014}.} and induces a negligible cost in real calculations. Since the maximal number of points is fixed as $s=1000$ in the octree, the S2T cost scales as $\Ocal(sN_B)$, where $N_B$ is the number of points in octree leaf boxes adjacent to the periodic boundaries in $B_0$. For the randomly distributed points in Figure~\ref{fig:points} B, $N_B\ll N$ and adding periodic boundary condition has little cost in $\tau_{\Ncal}$ as shown by Table~\ref{tab:StokesTimingRandom}. By reusing the octree data of $B_0$, the triply periodic cases only cost about $10\%\sim40\%$ more time than the free-space cases even in high precision cases, much more efficient than the previously proposed method \cite{gumerov_method_2014} where the near-field $\Ncal^{B_0}$ is calculated without reusing the octree for $B_0$.

The cost of calculating the far field
$\tau_{\Fcal}$ scales as $\Ocal(Tp^2)$, and we implemented step d) in
Algorithm~\ref{alg2} as a simple double loop over all target points
$\bx_{t}$ and downward equivalent sources $\by_{s}^{B_0,d}$, leaving
the optimization work to the \verb|C++| compiler. If the double loop
is further optimized with hand-tuned SIMD instructions including SSE,
AVX and FMA instructions, further speed up is possible. However, since
$\tau_{\Fcal}$ costs only about 5\% of the total computation time,
reduction in $\tau_{\Fcal}$ gives little benefit.

\subsection{Effect of $\ell$ }
\label{subsec:effectl}
In previous subsections, we fixed $\ell=2$. This is arbitrarily chosen because we found that changing $\ell$ has little effect on both timing and accuracy. In this subsection, we provide a comparison of different $\ell$ to illustrate this. All timing results here and in the following sections are collected on the same machine with the same settings as used in Section~\ref{subsec:stokestiming}.

If $\ell\gg 1$, the precomputation time $\tau_{M2L}^{pre}$ and $\tau_{\Ncal}^{pre}$ may both increase as discussed in the next section, which is an unnecessary cost. Further, in some cases if $\ell\gg 1$, the error may increase, because both $q_{\Ncal}^t$ and $q_{\Ncal}^t$ may be much larger than 1, and when added together $q_{\Ncal}^t+q_{\Ncal}^t$ may generate a large truncation error due to the finite floating point precision. Hence, we consider only $\ell\sim\Ocal(1)$.

\begin{figure}[h] 
	\centering
	\includegraphics[width=0.45\linewidth]{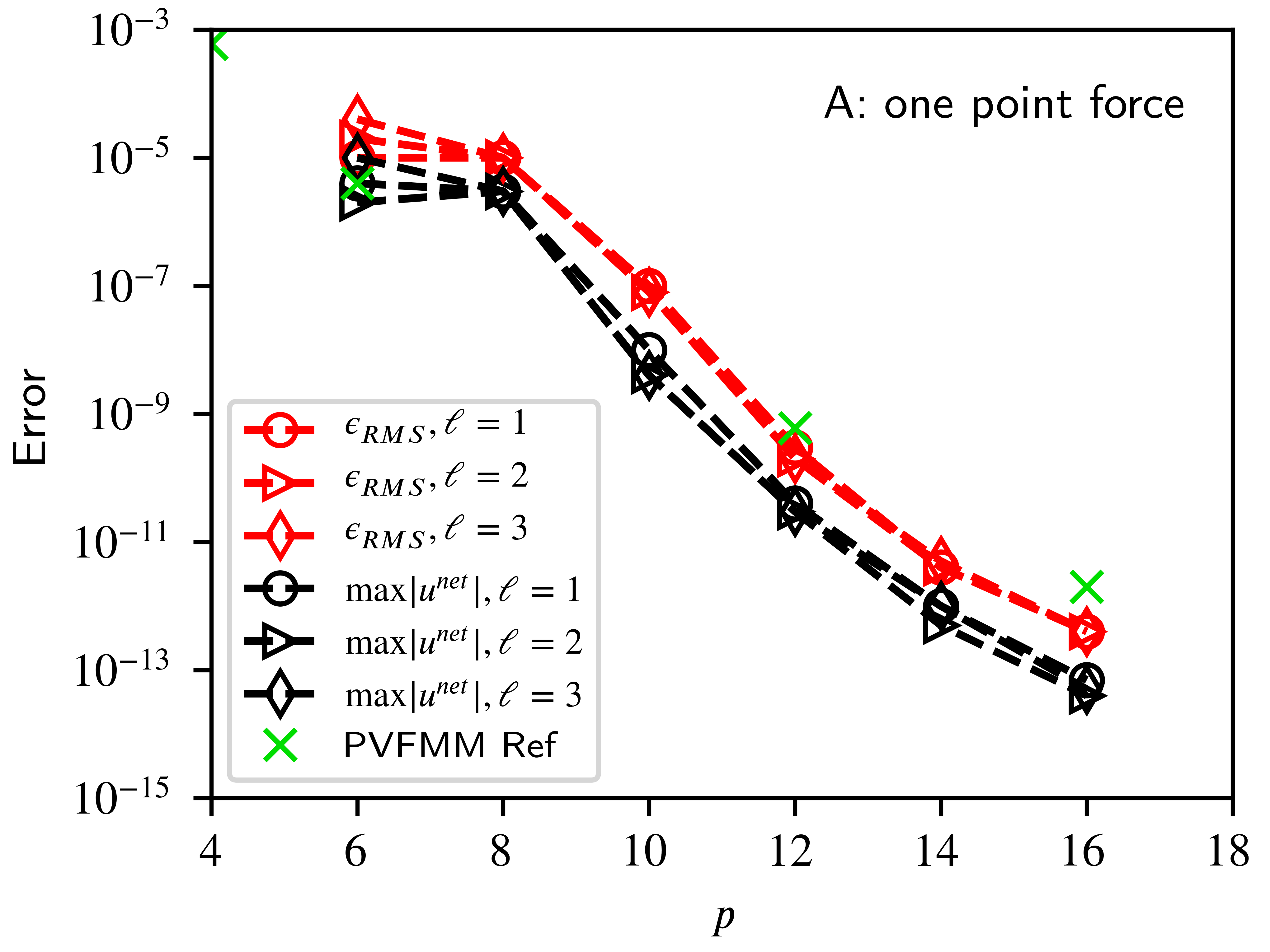}
	\includegraphics[width=0.45\linewidth]{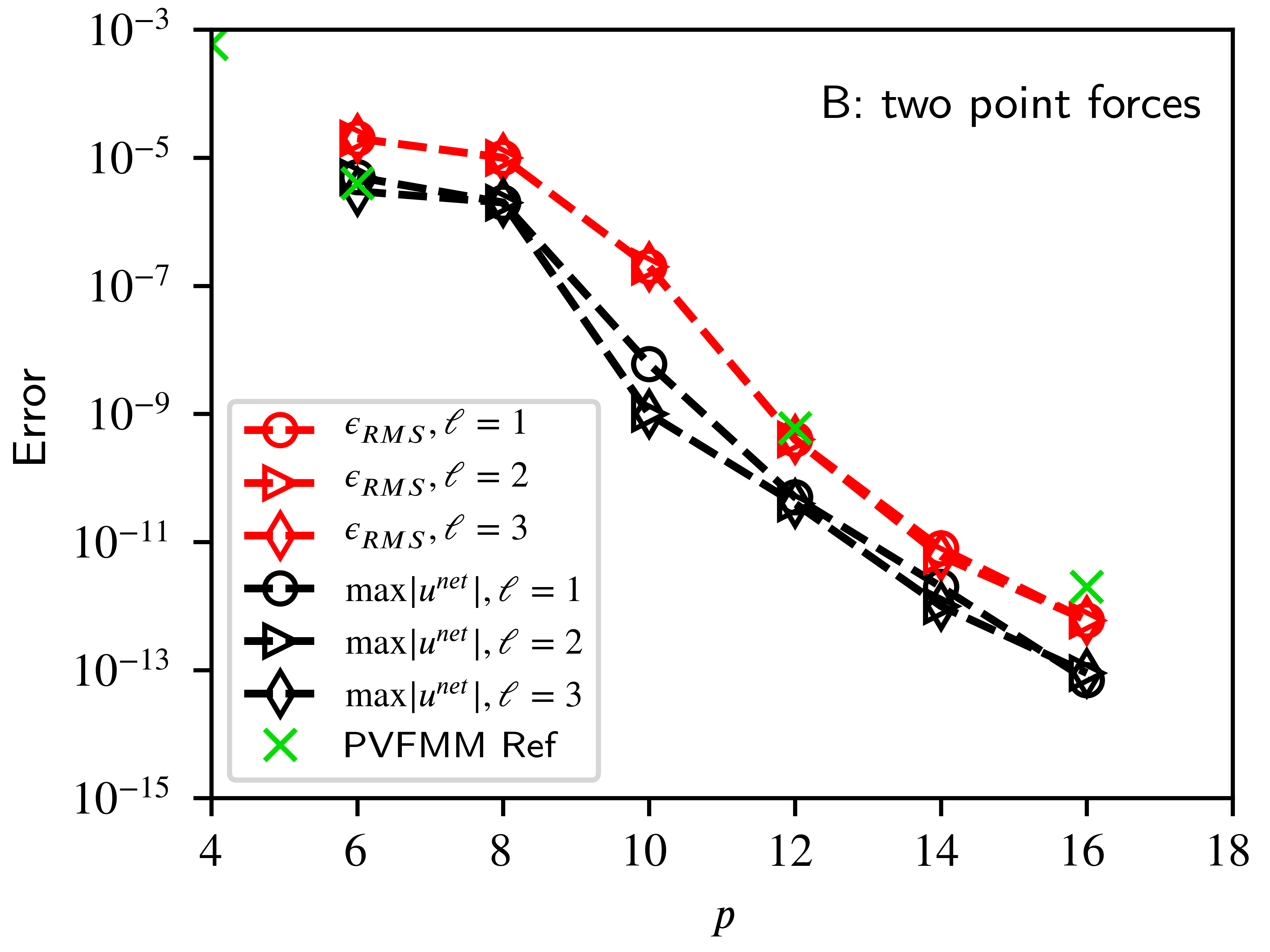}
	\caption{Accuracy test results for the triply periodic boundary condition with different $\ell$. The red symbols are point-wise	accuracy results, the black symbols are the no-flux condition accuracy results, and the green symbols are the reference error achieved by
		PVFMM \cite{malhotra_pvfmm_2015}. $\max\ABS{u^{net}}$ is the maximal component of $\bu^{net}$ in each case.}
	\label{fig:force3PEffectL}
\end{figure}

First, we repeated the accuracy test in Section~\ref{subsec:stokes3P} for $\ell=1,2,3$. The results are presented in Figure~\ref{fig:force3PEffectL}. It is clear that the accuracy results at different $\ell$ are almost undistinguishable, and they all satisfy the no-flux compatibility condition. 

\begin{figure} [h]
	\centering
	\includegraphics[width=0.45\linewidth]{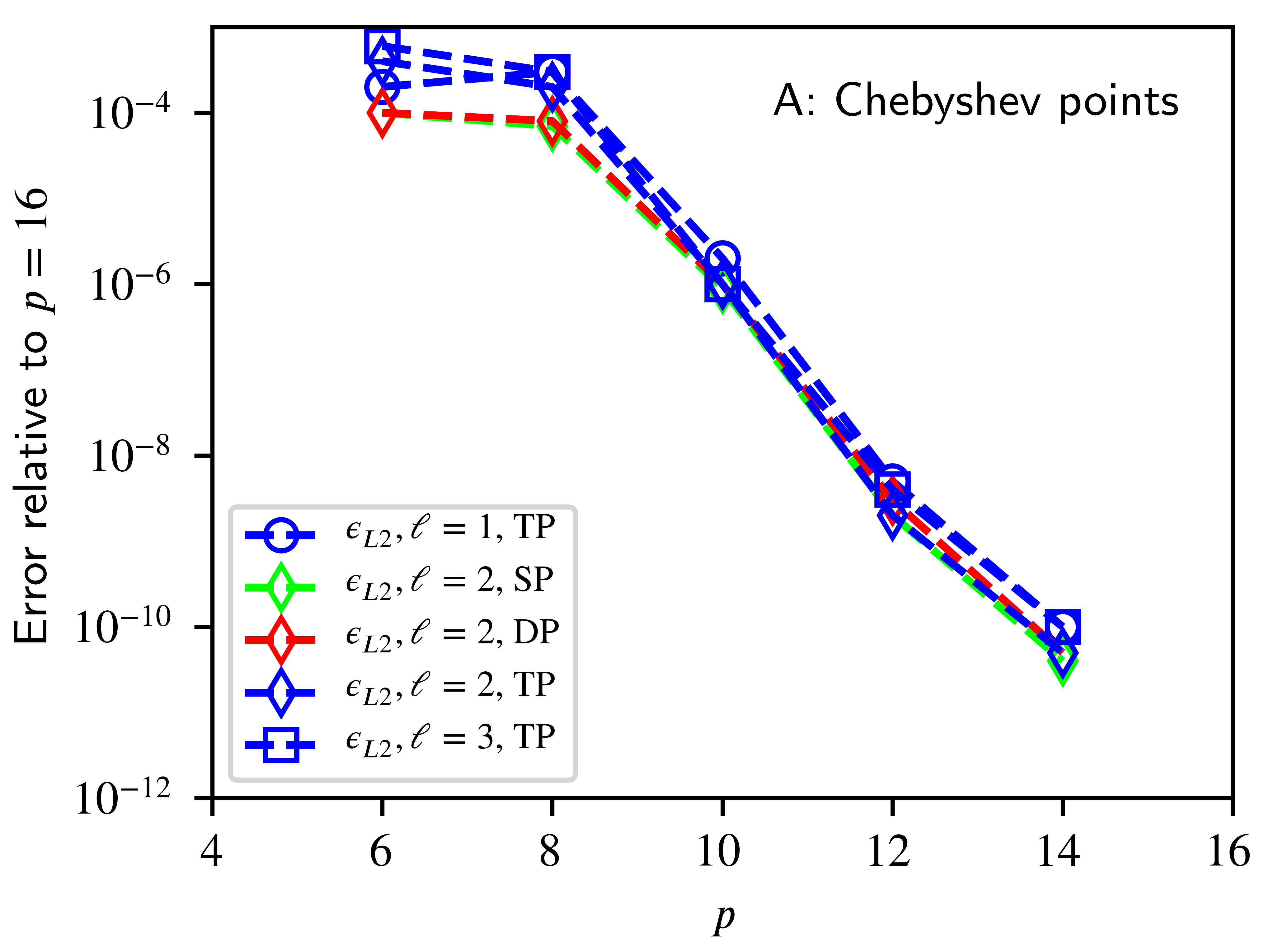}
	\includegraphics[width=0.45\linewidth]{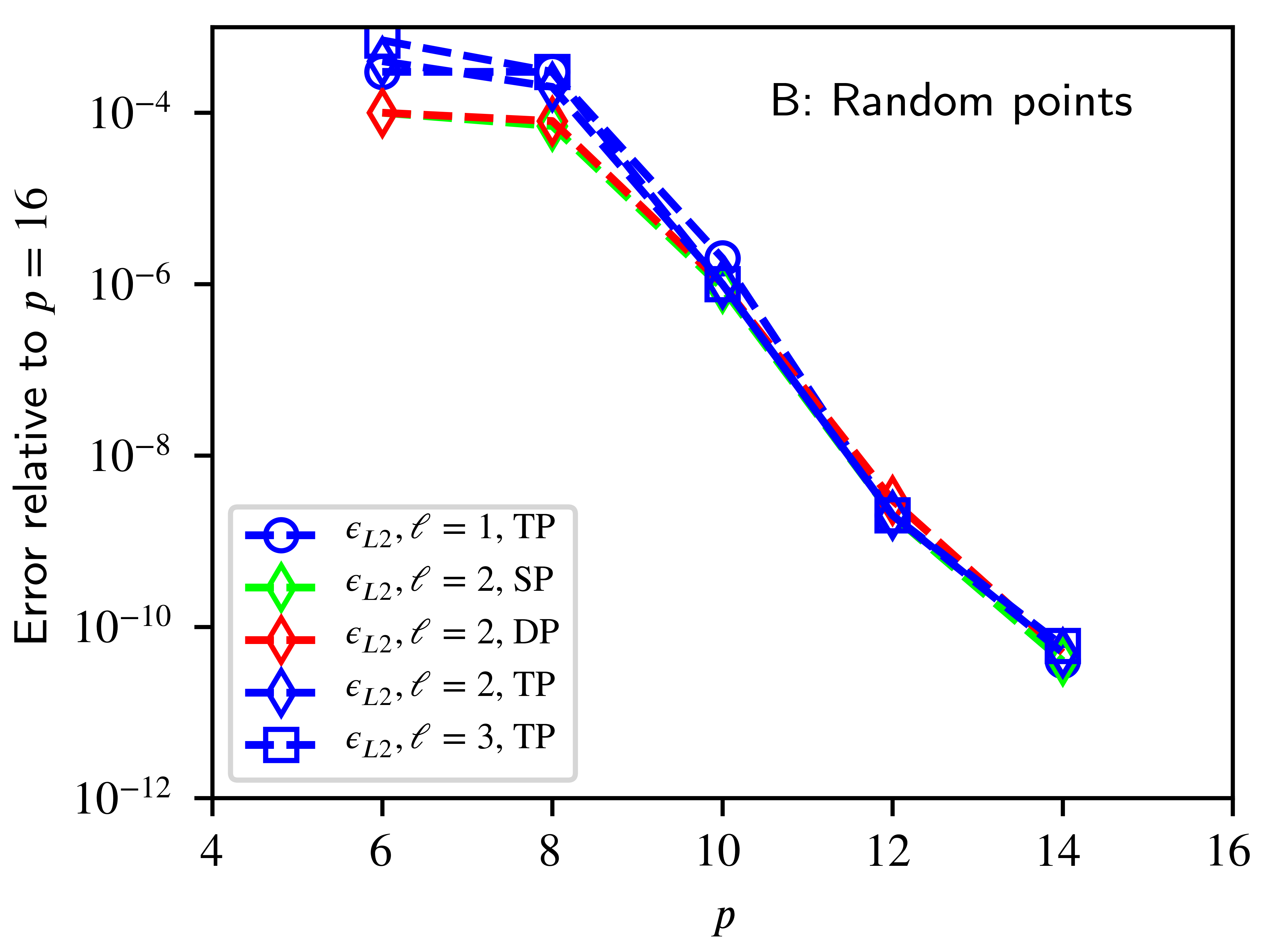}
	\caption{Convergence test results for different boundary condition with different $\ell$. The source and target points are the same set of $97^3$ Chebyshev and random points shown in Figure~\ref{fig:points} and used in the timing tests in Table~\ref{tab:StokesTiming} and~\ref{tab:StokesTimingRandom}. The red symbols are error compared to the results given by $p=16$.}
	\label{fig:ConvergeEffectL}
\end{figure}

\begin{table} [h]
	\centering
	\caption{Timing (in seconds) results for $\sim10^6$ Chebyshev
		points in Figure~\ref{fig:points} A for different $\Ncal$ sizes with the triply periodic boundary condition and Stokes velocity kernel. {As before, $\tau_{M2L}\sim 10^{-4}$ seconds and is much less than other parts, and so is not included in this table.} The tree construction time remains unchanged from the data shown in Table~\ref{tab:StokesTiming} and is also not included here.}
	\begin{tabular}{c|cc|cc|cc|c}
		\hline
		& \multicolumn{2}{c|}{ {$\ell=1$}} & \multicolumn{2}{c|}{ {$\ell=2$}} & \multicolumn{2}{c|}{ {$\ell=3$}} & HR \\
		$p$ & $\tau_{\Ncal}$ & $\tau_{\Fcal}$  & $\tau_{\Ncal}$ & $\tau_{\Fcal}$ &   $\tau_{\Ncal}$ & $\tau_{\Fcal}$ & \\
		\hline
		$6$ & 2.82 & 0.19 &2.92  &0.14 &   2.78 & 0.19 & 2.56 \\
		$8$ & 4.25 & 0.28  & 4.32 &0.25   &  4.36 & 0.27  & 4.16 \\
		$10$ & 5.72 & 0.43 & 5.67 &0.42   &   6.08 & 0.51 & 5.49 \\
		$12$ & 9.65 & 0.77 & 9.76  &0.70   &  9.72 & 0.77 & 9.39 \\
		$14$ & 16.85 & 1.00 & 16.67  &0.96  &  16.90 & 0.95 & 16.50 \\
		$16$ & 33.82 & 1.23 & 32.60  &1.24   &  32.70 & 1.23 & 32.34 \\
		\hline
	\end{tabular}
	\label{tab:StokesTimingEffectL}
\end{table}

Second, we repeated the two large tests for $\sim 10^6$ Chebyshev and random points in Section~\ref{subsec:stokestiming} to compare the convergence and timing for different $\ell$. The convergence error is calculated against the results with $p=16$ for each setting, and shown in Figure~\ref{fig:ConvergeEffectL}. Changing $\ell$ or changing periodicity has no discernible effects on the convergence error.

The timing results for different $\ell$ is shown in Table~\ref{tab:StokesTimingEffectL}, and we also added the timing for the Hierarchical Repetition (HR) method, which effectively uses $\ell=2^{30}$, from the original PVFMM package for comparison. It is worth noting that the results given by HR are off by a constant and do not satisfy the no-flux compatibility condition due to an unfixed constant flow across the root box $B_0$. For different $\ell$, $\tau_{\Fcal}$ remains almost constant as can be easily understood with Algorithm~\ref{alg2}. Further, $\tau_{\Ncal}$ is also almost constant and independent of $\ell$, and is also equal to the HR method. This is due to the way we precompute the KIFMM data and `compress' the near field evaluations into M2L translational operators, which originates from the HR method. This precompuation technique keeps $\tau_{\Ncal}$ independent of $\ell$, even for $\ell=2^{30}$ in the HR method. We shall discuss more about this precomputation in the next subsection.

\subsection{The precomputation of KIFMM and $\MLOP$}

\label{subsec:precompute}
\begin{figure}[h]
\centering 
\includegraphics[width=0.75\linewidth]{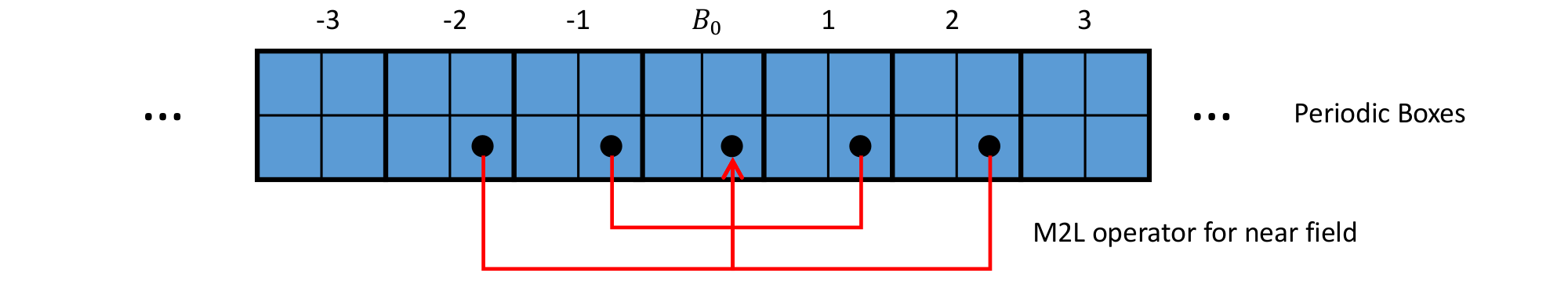}
\caption{A graphical representation for the precompute procedure for $\Ncal^{B_0}$. A translational M2L operator is directly summed for leaf boxes within $\ell$ layers of periodic images. The case for $\ell=2$ is shown in the figure. This is the same idea as the precomputation for the Hierarchical Repetition method, and more details can be found in the description of the PVFMM package \cite{malhotra_pvfmm_2015,malhotra_algorithm_2016}.}
\label{fig:precompute}
\end{figure}

With Algorithm~\ref{alg2}, the precomputation is split into two parts: the KIFMM evaluation and determining the $\MLOP$ operator. The former is for calculating the near field and the latter is for the far field. The precomputation for KIFMM is thoroughly optimized in the original PVFMM package \cite{malhotra_pvfmm_2015,malhotra_algorithm_2016} for the Hierarchical Repetition method. The method directly and hierarchically sum the multipole of each neighbor image box to get a M2L translational operator, and store it for later use. In real calculations, the KIFMM package only needs to read the operator, apply the M2L operator to the multipoles, and then travel down the tree. In this work, we replaced the hierarchical part of the summation with a simple loop of summation over the image boxes within $\ell$ layers of image boxes, as illustrated in Figure~\ref{fig:precompute}. \footnote{This modification is proposed by Wen Yan and developed by Dr. Malhotra, the original author of the PVFMM package.} After the modification the precomputation for applying boundary conditions takes much less time compared to the HR method since $\ell\sim \Ocal(1)$. 

The precomputation for $\MLOP$ is formulated in Section~\ref{sec:method}. In this work the precomputation time is not important so we used a pairwise Ewald method to calculate $\bK^{P,F}$, and a custom version of SVD. If the precompuation time for $\MLOP$ is important, accelerated methods such as the Spectral Ewald method \cite{lindbo_spectrally_2010} can be used.

The precomputation time for the Laplace kernel with double periodicity is shown in Table~\ref{tab:LaplaceDPTimingPrecomp}. $\tau_{\Ncal}^{pre}$ is the time needed by PVFMM itself and $\tau_{M2L}^{pre}$ is the time to construct $\MLOP$. For double periodicity the Ewald sum contains no hard-to-evaluate special functions so we used the Ewald sum to construct $K^{P,F}$. This greatly reduced the time compared to the direct summation method used for single periodicity in Section~\ref{subsec:laplace1D}. In general, $\tau_{\Ncal}^{pre}\approx \tau_{M2L}^{pre}$. Also $\tau_{\Ncal}^{pre}$ slightly increases with larger $\ell$, because to sum the M2L operator for near field the number of image boxes grows quadratically with increasing $\ell$ for double periodicity. It only slightly increases because the precomputation for boundary conditions and the M2L operator is not the only nor the dominating component, and a detailed comparison is beyond the scope of this work.

\begin{table}[h]
	\centering
	\caption{Comparison of precomputation time (in seconds) for Laplace kernel with double periodicity with different $\ell$.}
	\begin{tabular}{c|cc|cc|cc }
		\hline
		& \multicolumn{2}{c|}{ { $\ell=1$} } & \multicolumn{2}{c|}{ { $\ell=2$} } & \multicolumn{2}{c}{ { $\ell=3$} }  \\
		$p$  & $\tau_{\Ncal}^{pre}$ & $\tau_{M2L}^{pre}$  & $\tau_{\Ncal}^{pre}$ & $\tau_{M2L}^{pre}$  & $\tau_{\Ncal}^{pre}$ & $\tau_{M2L}^{pre}$  \\
		\hline
		$6$ & 0.17 & 0.25 &0.21 & 0.25 &0.23 &0.24 \\
		$8$ & 0.51 & 0.59 &0.53 & 0.59 &0.57 &0.59 \\
		$10$ & 1.11 & 1.27 &1.18 & 1.24 &1.26 &1.25 \\
		$12$ & 2.20 & 2.35 &2.31 & 2.37 &2.54 &2.31 \\
		$14$ & 4.15 & 4.21 &4.34 & 4.23 &4.80 &4.25 \\
		$16$ & 7.64 & 7.21 &8.83 & 7.30 &9.62 &7.31 \\
		\hline
	\end{tabular}
	\label{tab:LaplaceDPTimingPrecomp}
\end{table}

Since KIFMM is in general known to be slow for the Laplace kernel in comparison to the classical FMM, KIFMM is usually applied with other kernels. Here we report, in Table~\ref{tab:StokesTimingPrecomp}, the precomputation timing for the Stokes velocity kernel with triple periodicity. $\tau_{M2L}^{pre}$ is dominated by the SVD of matrix $\Acal$, followed by a dense matrix-matrix multiplication (DGEMM), and the cost for both parts scale as $\Ocal([\dim\Acal]^3)$. Since $\dim\Acal\propto p^2$, we have roughly $\tau_{M2L}^{pre} \propto p^6$. In real calculations, since the multi-threading overhead is relatively lower for larger matrices, the data shown in Table~\ref{tab:StokesTimingPrecomp} actually scales more close to $p^4$. This scaling is also true for $\tau_{\Ncal}^{pre}$ where SVDs and DGEMMs are also invoked. The precomputation time is compared to the HR method with effectively $\ell=2^{30}$, where no $\MLOP$ is constructed. In general, the hierarchical summation takes longer, and $\%_{HR} = (\tau_{\Ncal}^{pre} + \tau_{M2L}^{pre})/\tau_{HR}^{pre}\approx 60\%$. Besides satisfying the no-flux condition, the $\MLOP$ method also takes less time to setup, compared to the HR method. Further, if the precomputation time is important to applications, the pairwise Ewald can be replaced by Spectral Ewald method to further reduce $\tau_{M2L}^{pre}$.

The memory cost to construct $\MLOP$ is also dominated by the SVD of $\Acal$, which is a dense matrix with dimension shown in Table~\ref{tab:cond}. The memory cost scales as $p^4$ since $\dim\Acal\propto p^2$. For example, for Stokes velocity kernel with $p=16$, $\MLOP$ takes about 1GB of memory to compute. The details for this memory cost depends on the particular implementation and algorithm of the SVD procedure.

\begin{table}[h]
	\centering
	\caption{Comparison of precomputation time (in seconds) between the $\MLOP$ method and the Hierarchical Repetition (HR) method for the Stokes velocity kernel with triple periodicity. In HR, the hierarchy is repeated 30 times, generating a summation equivalent to $[2^{30}]^3$ images of $B_0$. $\%_{HR} = (\tau_{\Ncal}^{pre} + \tau_{M2L}^{pre})/\tau_{HR}^{pre}$ is the ratio of the total precomputation cost of $\MLOP$ method compared to the HR method.}
	\begin{tabular}{c|ccc|ccc|ccc|c}
		\hline
		& \multicolumn{3}{c|}{ { $\ell=1$} } & \multicolumn{3}{c|}{ { $\ell=2$} } & \multicolumn{3}{c|}{ { $\ell=3$} } & HR   \\
		$p$  & $\tau_{\Ncal}^{pre}$ & $\tau_{M2L}^{pre}$ & $\%_{HR}$ & $\tau_{\Ncal}^{pre}$ & $\tau_{M2L}^{pre}$ & $\%_{HR}$ & $\tau_{\Ncal}^{pre}$ & $\tau_{M2L}^{pre}$ & $\%_{HR}$ & $\tau_{HR}^{pre}$   \\
		\hline
		$6$ & 1.13 & 1.64 & 47\% & 1.18 & 1.69 & 49\% & 1.25 & 1.77 & 52\% & 5.8 \\
		$8$ & 3.32 &5.59 & 61\% & 3.43 & 5.59 & 62\% & 3.74 & 5.65 & 64\% & 14.63 \\
		$10$ & 10.33 &13.35 & 63\% & 11.69 & 13.51 & 66\% & 12.45 & 13.79 & 69\% & 38.00 \\
		$12$ & 26.89 &28.87 & 45\% & 29.19 & 29.68 & 48\% & 33.37 & 29.99 & 51\% & 123.92 \\
		$14$ & 69.42 &58.25 & 49\% & 73.30 &59.64 & 51\% & 83.16 & 60.45 & 55\% & 259.10 \\
		$16$ & 156.21 &110.04 & 57\% & 165.19 &110.61  & 59\% & 181.65 & 114.06 & 64\% & 465.45 \\
		\hline
	\end{tabular}
	\label{tab:StokesTimingPrecomp}
\end{table}

\section{Discussions and Conclusions}
\label{sec:discuss}
In this work we proposed a method to impose periodic boundary
conditions for KIFMM at a small cost. It can be implemented as an
efficient add-on layer for any KIFMM code, as long as $\Ncal^{B_0}$
within an \emph{integer} $\ell$ layers of neighboring periodic image
boxes can be calculated. The add-on periodizing layer adds the
contribution from $\Fcal^{B_0}$ to the results given by FMM routines,
and is straightforward to implement as a simple double loop at a small
cost of $\Ocal({Tp^2})$ (cf. Table~\ref{tab:StokesTiming}
and~~\ref{tab:StokesTimingRandom}), where $T$ is the number of target
points in $B_0$. On distributed memory clusters, the add-on layer
requires no inter-node communications and thus poses no difficulties
in scalability. The method is described in Algorithm~\ref{alg2} and is
demonstrated to require a small extra cost in addition to free-space
KIFMM. Further optimization is possible because with some
modifications of the underlying KIFMM routine, the evaluation of
$q_{\Fcal}^t$ can be combined with the KIFMM itself. In fact, $\MLOP$
can be incorporated into the L2L steps in the downward pass of the
octree data structure, and can be precalculated and stored as one
usually does for the KIFMM octree data precalculation step. In this
way, $\tau_{\Fcal}$ can be eliminated. However, since $\tau_{\Fcal}$
is already small, as shown in Table~\ref{tab:StokesTiming}
and~\ref{tab:StokesTimingRandom}, there is not much benefit of
implementing further optimizations.

In comparison to the hierarchical repetition method, the added value
of the $\MLOP$ method is properly handling of the compatibility
condition Eq.~(\ref{eq:compatcondition}). When the compatibility condition
is simple as in the example of 1D Madelung constant in
Section~\ref{subsec:laplace1D}, the hierarchical repetition method
should work equally well. However for more flexible compatibility
conditions the $\MLOP$ method is able to handle it properly as shown
in Section~\ref{subsec:stokes3P}. Another example of a non-trivial
compatibility condition is the well-known one-component plasma system,
where positive point charges move in a background of uniform rigid
neutralizing negative charge. In this case the neutrality
condition Eq.~(\ref{eq:compatcondition}) is not simply imposed on point
charges, but includes the background also. For one-component-plasma,
applying Ewald summation in $\MLOP$ gives a simple way to circumvent the
conditional convergence of the periodic charge summation. Also in this work we presented examples where the periodic
kernel $K^P$ is explicitly known. If not known, such periodic kernels
can be constructed in the precomputing stage for $\MLOP$ by inserting
control points \cite{barnett_unified_2016,marple_fast_2016}. Further,
this $\MLOP$ method works not only for KIFMM.
Our method works by adding a far field to the easy-to-calculate near field. It does not matter how the near field ($B_0$ with its image within $\ell$ layers) is calculated. 
If the near field is calculated by KIFMM, we can reuse the multipoles for the root box $\phi^{B_0,u}$ to construct the far field. 
Other methods suitable for the problem size, such as brute force summation or the Barnes-Hut tree can also be used. 
After all, the near field calculation is only a finite system with the open boundary condition, and is straightforward to calculate. With those alternative methods other than KIFMM, we may need to calculate $\phi^{B_0,u}$ separately at some extra cost because there is no data to reuse.

The method requires a precomputation of $\MLOP$, which is found from the
free-space kernel $\bK$, the `periodic' kernel $\bK^{P}$, and the
known discretization and locations of check and equivalent
surfaces. For typical ``$1/r$'' kernels $\bK$, such as the Laplace or
Stokes kernels, $\bK^{P}$ is usually expressed as the Ewald form
$\bK^E$, and has been thoroughly discussed in the literature for
singly, doubly and triply periodic systems. In fact, $\MLOP$ can be
calculated for different far-field geometries, and so is not limited
to periodic systems. For example, the far-field can be calculated to
represent a channel with finite instead of infinite length, or in the
presence of a reflecting wall with the method of images. If the
requirement of accuracy can be relaxed, $\MLOP$ can be used to
approximate a ring-channel geometry to study the spontaneous
coherent flows created by active
suspensions \cite{wioland_directed_2016,giomi_one_2016}. The flow of liquids in fractal
geometries  \cite{heinen_classical_2015} is another example of
interesting problems that could be studied.

In this work we focused on a cubic domain $B_0$, and all shapes which
fit into a cubic box $B_0$ can be treated without modification, for
example, a long tube with singe periodicity in the tube direction. In
principal, slab geometries, where the size of $B_0$ in a non-periodic
direction is much larger than that in the periodic direction, can be
handled with the same $\MLOP$ method. When the aspect ratio of $B_0$ is of order $\Ocal(1)$, Ewald methods can still be used. For higher but still moderate aspect ratios, Ewald methods may converge slowly but we can utilize the scheme in \cite{barnett_unified_2016} to construct the periodic kernel $\bK^P$ directly without a series summation, with some tuning of the discretization of the check and equivalent surfaces to fit the high aspect
ratio. For extremely high aspect ratio the accurate construction of the periodic kernel $\bK^P$ remains to be studied. Another important case is a cubic domain $B_0$ sheared along
one axis, which is widely used in the simulation of particle
suspension dynamics to retrieve rheology
information \cite{wang_spectral_2016}. However, to the best of our
knowledge, due to the complexity of building octrees there is no high
performance KIFMM code available to academia that can handle a box of
arbitrary aspect ratio and sheared deformation. Therefore, although
the idea of $\Ncal,\Fcal$ splitting and $\MLOP$ operator are still
useful, the underlying FMM remains lacking, and this is well worth
further study.

In recent years, many algorithms and applications are proposed for the
algebraic variants of
FMM \cite{minden_technique_2016,Coulier2015,ibeid_matrix-free_2016,sun_matrix_2001,yokota_fast_2016}. The
algebraic version utilizes the hierarchical low-rank structure of FMM
matrix $\bm{\Mcal}_{t,s}$: $\bq_t=\bm{\Mcal}_{t,s}\bphi_s$. The
low-rank structure comes from the fact that for a box $B_i$ in an
octree, potentials due to the far-field $\Fcal^{B_i}$ can be
approximated by equivalent sources with a small number of degree of
freedom, and therefore the interaction between these two leaf boxes
can be compressed by a low rank matrix. It is hierarchical in the
sense that for all levels of boxes in the octree, corresponding low
rank structures can be formed. With the $\MLOP$ method, the FMM matrix
is split into near and far field:
\begin{align}
	\bm{\Mcal}_{t,s} = \bm{\Mcal}_{t,s}^{\Ncal^{B_0}} + \bm{\Mcal}_{t,s}^{\Fcal^{B_0}}
\end{align}  
where $\bm{\Mcal}_{t,s}^{\Ncal^{B_0}}$ is simply a FMM matrix and
preserves the structure. $ \bm{\Mcal}_{t,s}^{\Fcal^{B_0}}$ is a
low-rank matrix since it is applied through $\phi^{B_0,d}$ and
$\bK^{P,F}$. Therefore, our method should be compatible with the
algebraic variants of FMM. Another possible extension is to combine
Algorithm~\ref{alg2} with the Quadrature-by-Expansion (QBX)
method \cite{klockner_quadrature_2013} to evaluate near-singular kernel
functions in boundary layer integral problems.

\section{Acknowledgement}
WY thanks Dhairya Malhotra for his help in modifying his PVFMM
package to implement $\Ncal^{B_0}$ calculation. MJS acknowledges support from NSF grants DMS-1463962 and DMS-1620331, and NIH Grant GM104976. WY and MJS thank Alexander Barnett and Shravan Veerapaneni for inspiring our idea of periodizing KIFMM, and for their critical readings of the manuscript. 

\appendix
\section{Open-sourced implementation}
An open-sourced implementation of the method described in this paper
will be available at the GitHub. The modified PVFMM to calculate
$q_{\Ncal}^t$ with SP, DP and TP boundary conditions will be available
at the author's fork of PVFMM. An implementation to calculate the
periodizing operator $\MLOP$ and a wrapper library to calculate full
periodic FMM for the Stokes velocity kernel will be available in a
separate repository. A matlab script based on Spectral-Ewald
method will also be available to help the readers quickly implement
$\MLOP$ for other Stokes and Laplace kernels with the Spectra-Ewald
method.

\section{The backward stable solver}
\label{sec:svdcompare}

Table~\ref{tab:svd} compares various implementations of the backward
stability solver, where $\Acal$ is a square matrix of dimension
$4056\times4056$ for the Stokes velocity kernel with $p=16$. It is the
most challenging case for all $\Acal$ in Table~\ref{tab:cond}. We
calculate the max absolute backward error of solving for $\bx^*$
satisfying $\Acal \bx = \bone$: $\max\epsilon_{abs} =
\max\ABS{\Acal\bx^*-\bone}$, where $\bone$ is a vector with every
entry being $1$. SVDpvfmm is the SVD solver implemented in
PVFMM \cite{malhotra_pvfmm_2015}. QR is the
\verb|fullPivHouseholderQr().solve()| routine in the popular C++
linear algebra library \verb|Eigen 3.3.3|. JacobiSVD is the
\verb|JacobiSVD| routine with
\verb|FullPivHouseholderQRPreconditioner| in
\verb|Eigen|. $\Acal^+$ is the same as JacobiSVD, but
$\Acal^{\dagger}$ is explicitly calculated. DGESDD uses the
\verb|lapacke_dgesdd| routine from \verb|LAPACK 3.7|. DGESDD0 means
the threshold $\epsilon_{SVD}=0$, and in other routines
$\epsilon_{SVD}$ is set to $\epsilon_{m}\max \ABS{s}$, where
$\epsilon_{m}$ is the machine accuracy and $\max \ABS{s}$ is the max
absolute value of all singular values.  A Jacobi SVD routine
\verb|lapacke_dgejsv| is also available from \verb|LAPACK|, which is
contributed by the inventor of the
algorithm \cite{drmac_new_2008,drmac_new_2008_2}.

\begin{table}[ht]
	\centering
	\caption{Survey of SVD for solving $\Acal\bx=\bone$ for the Stokes velocity kernel with $p=16$. } 

	\begin{tabular}{ccccccc}
		\hline
		       &  SVDpvfmm  & QR  & JacobiSVD & $\Acal^{+}$  & DGESDD & DGESDD0 \\
		\hline
		$\max\epsilon_{abs}$ & $1.1E-13$   & $2.5E-13$ & $7.6E-14$ & $2.3E-09$ & $3.0E-14$ & $5.5E-14$ \\
		\hline
	\end{tabular}
	\label{tab:svd}
\end{table}

The results in Table~\ref{tab:svd} show that backward stability can be
achieved with various SVDs with $\bU,\bV$ stored separately. However,
we found that the speed of those solvers vary significantly. SVDpvfmm
and DGESDD are partially threaded with \verb|OpenMP| and are the
fastest. JacobiSVD from \verb|Eigen| is significantly slower and is
not threaded, but is still acceptable because we need to calculate the
SVD once for all columns in $\MLOP$. In the results presented in this
work we used SVDpvfmm to calculate $\MLOP$, and it is demonstrated in
the accuracy tests that we achieved the bound of accuracy $10^{-13}$.

\bibliographystyle{elsarticle-num}
\section*{\refname}

\end{document}